\documentclass[10pt]{article}

\usepackage{amsmath,amsfonts,graphics,graphicx,bm,enumerate,epsfig,psfrag,subfig}
\usepackage{palatino,cite}
\usepackage{color,soul}
\usepackage{xcolor}
\usepackage[mathscr]{euscript} 
\usepackage{mathabx}

\usepackage{float}
\usepackage[colorlinks=true, pdfstartview=FitV, linkcolor=black, citecolor= black, urlcolor= black]{hyperref}
\usepackage{footnpag}			       

\usepackage{siunitx}
\usepackage{longtable,tabularx}

\newcommand{\m}[1]{{\bf{#1}}}
\newcommand{\g}[1]{\boldsymbol #1}

\newcommand{\tr}{^{\sf T}}
\newcommand{\C}[1]{{\cal {#1}}}

\newcommand{\ddtm}[2]{\hspace{-1mm} ^{^{^{^ #1}}} \hspace{-1.5mm} \frac{d}{dt} \left( {#2} \right)}
\newcommand{\om}[2]{{^\C{#1} \hspace{-0.1mm} \g{\omega} \hspace{-0.1mm} ^\C{#2}}}
\newcommand{\al}[2]{{^\C{#1} \hspace{-0.1mm} \g{\alpha} \hspace{-0.1mm} ^\C{#2}}}
\newcommand{\sk}{^{\textrm{x}}}
\newcommand{\br}[2]{{{\left\lbrace {#2} \right\rbrace}_{#1}}}
\newcommand{\brsk}[2]{{{\left\lbrace {#2} \right\rbrace}_{#1} \sk}}
\newcommand{\brtr}[2]{{{\left\lbrace {#2} \right\rbrace}_{#1} \tr}}
\newcommand{\w}[1]{\omega_{#1}}
\newcommand{\e}[1]{\epsilon_{#1}}
\newcommand{\edot}[1]{\dot{\epsilon}_{#1}}
\newcommand{\h}[1]{{\m{#1}}}

\newcommand{\lsup}[2]{{\vphantom{#2}}^{#1}{{#2}}}

\newcommand{\CBA}{\m{C}_{\C{B} \C{A}}}
\newcommand{\CAB}{\m{C}_{\C{A} \C{B}}}

\newcommand{\CAE}{\m{C}_{\C{A} \C{E}}}
\newcommand{\CEA}{\m{C}_{\C{E} \C{A}}}

\newcommand{\CBE}{\m{C}_{\C{B} \C{E}}}
\newcommand{\CEB}{\m{C}_{\C{E} \C{B}}}

\newcommand{\singlespacing}{\renewcommand{\baselinestretch}{1}\normalsize\normalfont}

\date{}
\begin{document}

\title{\bf Nonsingular Parameterization for Modeling \\ Translational Motion Using Euler Parameters}

\author{Alexander T.~Miller\thanks{Ph.D.~Student, NDSEG Fellow, Department of Mechanical and Aerospace Engineering.  E-mail:  alexandertmiller@ufl.edu.} \\ Anil~V.~Rao\thanks{Professor, Erich Farber Faculty Fellow, and University Term Professor, Department of Mechanical and Aerospace Engineering.  E-mail:  anilvrao@ufl.edu.  Associate Fellow AIAA.  Corresponding Author.} \\\\ {\em University of  Florida} \\ {\em Gainesville, FL 32611}}

\maketitle{}
\singlespacing

 \begin{abstract}
   A parameterization is described for quantifying translational motion of a point in three-dimensional Euclidean space.  The parameterization is similar to well-known parameterizations such as spherical coordinates in that both position and velocity are decoupled into magnitude and orientation components.  Unlike these standard parameterizations, where principal rotation sequences are employed, the method presented in this research employs Euler parameters.  By using Euler parameters instead of Euler angles, singularities and trigonometric functions are removed from the equations of motion. The parameterization is demonstrated on two examples, where it is found that the new parameterization offers both mathematical and computational advantages over other commonly used parameterizations.  
\end{abstract}

\section{Introduction}
Describing the motion of particles and rigid bodies moving in three-dimensional Euclidean space is fundamental to the engineering, physics, and mathematics communities.  The subject of particle and rigid body mechanics is divided into two distinct but related parts:  kinematics and kinetics.  In particular, kinematics is a critical part of developing the equations that describe motion of a mechanical system, and a key aspect of kinematics is the choice of coordinates used to parameterize the motion.  Classical kinematic parameterizations of particle motion include Cartesian, cylindrical, and spherical coordinates.


In recent years, renewed attention has been given toward developing coordinate systems to describe the translational motion of a point in three-dimensional Euclidean space.  A great deal of this research has focused on the use of quaternions.  While traditionally unit quaternions (also known as Euler parameters) are used to describe the orientation of a rigid body \cite{Shuster1993, Kane1983, Wie1985, Wie2008, Hughes2, Junkins2012, Wertz2012}, more recent work has focused on extending quaternions to describe translational motion \cite{KS1965, Chelnokov1981, Vivarelli1983, Deprit1994, Vrbik1994, Vrbik1995, Waldvogel2006, Waldvogel2008, Saha2009, Broucke1971, Deprit1975, Gurfil2005, Chelnokov2001, Chelnokov2003, Chelnokov2013, Chelnokov2014, Chelnokov2019, Libraro2014, Roa2017}.  In particular, the seminal work of Ref.~\cite{KS1965} used spinors to model three-dimensional motion in four dimensions.  Although the work of Ref.~\cite{KS1965} recognized the similarity between their KS transformation and quaternion multiplication, it was not until the works of Refs.~\cite{Chelnokov1981, Vivarelli1983} that the connection was made between the KS transformation and quaternions.  The connection between the KS transformation and quaternions has since been revisited and expanded upon in the works of Refs.~\cite{Deprit1994, Vrbik1994, Vrbik1995, Waldvogel2006, Waldvogel2008, Saha2009}.  In addition to the KS transformation, quaternions have been introduced to other translational motion descriptions.  The first three angles in the set of classical orbital elements, longitude of ascending node, orbital inclination, and argument of periapsis, form a $3-1-3$ Euler angle sequence which is equivalently expressed by Euler parameters to eliminate the singularty associated with an equatorial orbit \cite{Broucke1971, Deprit1975, Gurfil2005}.  Finally, Refs.~\cite{Chelnokov2001, Chelnokov2003, Chelnokov2013, Chelnokov2014, Chelnokov2019,Libraro2014,Roa2017} have further developed parameterizations using quaternions for orbital motion.

This research builds upon the aforementioned works using quaternions to describe the translational motion of a point in three-dimensional Euclidean space.  In the prior work, position is typically decoupled into magnitude and orientation components by employing a single unit quaternion.  Different from the prior work, this research employs {\em two} sets of Euler parameters (two unit quaternions) such that {\em both} the position and the velocity of a point in three-dimensional Euclidean space are decoupled into separate variables for magnitude and orientation.  In particular, the magnitudes are given by the radial distance, $r$, and the speed, $v$.  Likewise, the Euler parameters define the position and velocity directions.  The complete set of parameters is termed here as "$rv$-Euler parameters".

Conceptually, the $rv$-Euler parameterization is similar to a spherical coordinate parameterization.  Both $rv$-Euler parameters and spherical coordinates include the radial distance $r$ and the speed $v$ as variables.  Furthermore, each of these parameterizations employs separate variables to define the directions of position and velocity.  Thus, velocity dependent forces (for example, lift and drag) are easily characterized in these formulations.  However, the $rv$-Euler parameters differ from spherical coordinates in that Euler parameters, not Euler angles, define the position and velocity directions.  In this regard, the $rv$-Euler parameterization is similar to the quaternionic forms of the equations of motion developed in Refs.~\cite{Chelnokov2001, Chelnokov2003, Libraro2014, Roa2017}.  Employing Euler parameters instead of Euler angles yields the computational advantage that the equations of motion do not contain trigonometric functions or singularities (for example, the polar latitude singularities using spherical coordinates).  Thus, the $rv$-Euler parameterization combines the computational benefits afforded by Euler parameters with the added benefit of describing systems involving {\em both} position and velocity dependent forces.

The remainder of the paper is organized as follows.  Section~\ref{sect:NC} details the notation and conventions which will be employed throughout the paper.  A synopsis of the relevant math is provided in Section~\ref{sect:Review}.  The $rv$-Euler parameters are introduced in Section~\ref{sect:Derivation1} followed by the derivation of the equations of motion.  The new equations of motion are demonstrated on two examples in Section~\ref{sect:Examples}.  Finally, the key aspects of the $rv$-Euler parameterization are discussed in Section~\ref{sect:Discussion} and Section~\ref{sect:Conclusion} contains conclusions of the research.

\section{Notation and Conventions \label{sect:NC}}
The following notation and conventions will be employed throughout the paper.  Scalars will be represented by lowercase symbols (for example, $p \in \mathbb{R}$).  Next, vectors lie in three-dimensional Euclidean space and will be denoted by lowercase bold symbols (for example, $\m{p} \in \mathbb{E}^3$).  Moreover, the usual notation ``$\cdot$" and ``$\times$" specify scalar and vector products.  

Next, reference frames will be represented by uppercase calligraphic letters (for example, $\C{A}$).  Each reference frame will be designated a single coordinate system fixed in the frame.  Moreover, each coordinate system will be comprised of an origin and a set of three, right-handed, orthonormal basis vectors.  Table~\ref{tab:frames} summarizes the notation for the various reference frames, origins, and basis vectors which will be used.  The precise definitions of the basis vectors will be laid out in Section~\ref{sect:Derivation1}.

\begin{table}[h]
  \centering
  \caption{Reference Frames.\label{tab:frames}}
  \renewcommand{\baselinestretch}{1}\normalsize\normalfont
 \begin{tabular}{|c|c|c|c|}\hline
    \textbf{Frame} & \textbf{Notation} & \textbf{Basis} & \textbf{Origin} \\\hline
    Inertial             & $\C{N}$     & $\{ \h{n}_1, \h{n}_2, \h{n}_3 \}$ & $O$ \\
    Observation & $\C{E}$      & $\{ \h{e}_1, \h{e}_2, \h{e}_3 \}$   & $O$ \\
    Position           & $\C{A}$      & $\{ \h{a}_1, \h{a}_2, \h{a}_3 \}$   & $O$ \\
    Velocity           & $\C{B}$      & $\{ \h{b}_1, \h{b}_2, \h{b}_3 \}$   & $P$  \\\hline
  \end{tabular}
\end{table}

Next, consider the expressions of vectors in particular bases.  The notation $\br{\C{A}}{\m{p}}$ will be employed to denote the expression of a vector $\m{p}$ in the (unique) basis assigned to frame $\C{A}$.  Moreover, $\br{\C{A}}{\m{p}}$ is equivalent to the column matrix $[p_1~p_2~p_3]\tr \in \mathbb{R}^3$ if $\m{p} = p_1 \h{a}_1 + p_2 \h{a}_2 + p_3 \h{a}_3$.  Furthermore, all direction cosine matrices will be represented by an uppercase bold ``$\m{C}$" with a subscript denoting the relevant frames and indicating the direction of the transformation.  Note that the relevant basis vectors are implied by the frame as laid out in Table~\ref{tab:frames}.  Thus, the identities
\begin{equation}\label{eq:DCM-def}
\begin{array}{lclclcl}
\br{\C{B}}{\m{p}} &=& \CBA \br{\C{A}}{\m{p}} & \textrm{and} & \br{\C{A}}{\m{p}} &=& \CAB \br{\C{B}}{\m{p}},
\end{array}
\end{equation}
hold for any arbitrary vector $\m{p}$, choice of frames $\C{A}$ and $\C{B}$, and choice of right-handed, orthonormal bases $\{ \h{a}_1, \h{a}_2, \h{a}_3 \}$ for frame $\C{A}$ and $\{ \h{b}_1, \h{b}_2, \h{b}_3 \}$ for frame $\C{B}$.  It is noted that $\CAB = \CBA\tr = \CBA^{-1}$ as well.  Finally, the notations $\CAB(i,:)$, $\CAB(:,j)$, and $\CAB(i,j)$ will be used to denote, respectively, the $i^{th}$ row of $\CAB$, the $j^{th}$ column of $\CAB$, or element $(i,j)$ of $\CAB$.

Next, scalar and vector products of two vectors will be expressed in the basis of a desired frame as follows.  Let $\m{p}$ and $\m{q}$ be any two vectors and let $\C{A}$ be the desired frame.  Moreover, let $\br{\C{A}}{\m{p}} = [p_1~p_2~p_3]\tr$ and $\br{\C{A}}{\m{q}} = [q_1~q_2~q_3]\tr$ be the expressions of vectors $\m{p}$ and $\m{q}$ in the basis of frame $\C{A}$. Then, the scalar product of $\m{p}$ with $\m{q}$ expressed in the basis of frame $\C{A}$ will be denoted by $\br{\C{A}}{\m{p} \cdot \m{q}} = \brtr{\C{A}}{\m{p}} \br{\C{A}}{\m{q}}$.  Likewise, the vector product of $\m{p}$ with $\m{q}$ expressed in the basis of frame $\C{A}$ will be represented by $\br{\C{A}}{\m{p} \times \m{q}} = \brsk{\C{A}}{\m{p}} \br{\C{A}}{\m{q}}$, where
\begin{equation}\label{eq:skew}
\brsk{\C{A}}{\m{p}} = \begin{bmatrix}
0      & -p_3 & p_2\\
p_3   & 0      & -p_1\\
-p_2 & p_1   & 0
\end{bmatrix},
\end{equation}
denotes the skew-symmetric operator $\{ \cdot \}\sk$ acting on $\br{\C{A}}{\m{p}}$.  Finally, the identities
\begin{equation}\label{eq:dot-cross-express}
\begin{array}{lclclcl}
\brtr{\C{B}}{\m{p}} \br{\C{B}}{\m{q}} &=& \brtr{\C{A}}{\m{p}} \br{\C{A}}{\m{q}} & \textrm{and} & \brsk{\C{B}}{\m{p}} \br{\C{B}}{\m{q}}&=& \CBA \brsk{\C{A}}{\m{p}} \br{\C{A}}{\m{q}},
\end{array}
\end{equation}
relate the expressions of the scalar and vector products of vectors $\m{p}$ and $\m{q}$ in the bases of any two frames $\C{A}$ and $\C{B}$.

Lastly, the rate of change of scalars and vectors are defined as follows. First, all rates of change of scalars are denoted by the overdot symbol.  For example, the notation $\dot{p} = dp/dt$ denotes the rate of change of the scalar $p$.  Next, all rates of change of vectors are denoted with a left superscript to indicate the reference frame in which the rate of change is observed.  For example, the notation $\lsup{\C{A}}{d\m{p}/dt}$
denotes the rate of change of the vector $\m{p}$ as viewed by an observer in frame $\C{A}$.  Moreover, the rate of change of the vector $\m{p}$ as viewed by an observer in frame $\C{A}$ is related to the rate of change of the vector $\m{p}$ as viewed by an observer in frame $\C{B}$ via the transport theorem $\lsup{\C{A}}{d\m{p}/dt} = \lsup{\C{B}}{d\m{p}/dt} + \om{A}{B} \times \m{p}$,
where $\om{A}{B}$ is the angular velocity of frame $\C{B}$ as viewed by an observer in frame $\C{A}$.  Similarly, the angular acceleration of frame $\C{B}$ as viewed by an observer in frame $\C{A}$ is denoted by $\al{A}{B}$.  Finally, the velocity and acceleration as viewed by an observer in reference frame $\C{A}$ are denoted, respectively, as $\lsup{\C{A}}{\m{v}}$ and $\lsup{\C{A}}{\m{a}}$, where $\C{A}$ denotes the frame in which the rate of change is observed.  


\section{Mathematical Preliminaries \label{sect:Review}}
A brief review of Euler parameters is now presented that will be used later to develop the $rv$-Euler parameterization.
Consider a direction cosine matrix $\CBA$ where frames $\C{A}$ and $\C{B}$ (and their right-handed, orthonormal bases) are arbitrary.  It is known that a right-handed orthonormal basis $\{ \h{a}_1, \h{a}_2, \h{a}_3 \}$ fixed in frame $\C{A}$ can be rotated to be aligned with the basis $\{ \h{b}_1, \h{b}_2, \h{b}_3 \}$ fixed in frame $\C{B}$ by a single rotation about an axis $\h{q}$ by an angle $\phi$.  Suppose further that $\h{q}$ is a unit vector which, when expressed in the basis $\{ \h{b}_1, \h{b}_2, \h{b}_3 \}$, is given as $\br{\C{B}}{\h{q}} = [q_1, q_2, q_3]\tr$, noting that, for the specific case where $\m{q}$ is the axis of rotation, $\br{\C{B}}{\h{q}} = \br{\C{A}}{\h{q}}$.  Then, the axis-angle representation of $\CBA$ is written as
\begin{equation}\label{eq:C_BA-AxisAngle}
\CBA = \begin{bmatrix}
(1 - \cos \phi) q_1^2 + \cos \phi &
(1 - \cos \phi) q_1 q_2 + q_3 \sin \phi &
(1 - \cos \phi) q_1 q_3 - q_2 \sin \phi \\
(1 - \cos \phi) q_2 q_1 - q_3 \sin \phi &
(1 - \cos \phi) q_2^2 + \cos \phi &
(1 - \cos \phi) q_2 q_3 + q_1 \sin \phi \\
(1 - \cos \phi) q_3 q_1 + q_2 \sin \phi &
(1 - \cos \phi) q_3 q_2 - q_1 \sin \phi &
(1 - \cos \phi) q_3^2 + \cos \phi
\end{bmatrix}.
\end{equation}

The axis-angle parameters give rise to a physical interpretation of the Euler parameters which are defined from the axis-angle parameters as
\begin{equation}\label{eq:Quat}
\begin{array}{lclclclclclclcl}
\e{1} & = & q_1 \sin \phi/2 & , & \e{2} & = & q_2 \sin \phi/2 & , & \e{3} & = & q_3 \sin \phi/2 & , & \eta & = & \cos \phi/2.  
\end{array}
\end{equation}
In Eq.~\eqref{eq:Quat} it is observed that the Euler parameters $\{ \e{1}, \e{2}, \e{3} \}$ provide the same directional information as the axis-angle parameters $\{ q_1, q_2, q_3 \}$.  Note also that the Euler parameters satisfy $\e{1}^2 + \e{2}^2 + \e{3}^2 + \eta^2 = 1$
which is equivalent to a unit norm constraint on the quaternion defined by the vector part $[\e{1}, \e{2}, \e{3}]\tr$ and scalar part $\eta$.  Now, given the definitions for the Euler parameters in Eq.~\eqref{eq:Quat}, the description of $\CBA$ in Eq.~\eqref{eq:C_BA-AxisAngle} becomes
\begin{equation}\label{eq:C_BA-Quat}
\CBA = \begin{bmatrix}
1 - 2(\e{2}^2 + \e{3}^2) &
2(\e{1}\e{2} + \e{3} \eta) &
2(\e{1}\e{3} - \e{2} \eta) \\
2(\e{2}\e{1} - \e{3} \eta) &
1 - 2(\e{3}^2 + \e{1}^2) &
2(\e{2}\e{3} + \e{1} \eta) \\
2(\e{3}\e{1} + \e{2} \eta) &
2(\e{3}\e{2} - \e{1} \eta) &
1 - 2(\e{1}^2 + \e{2}^2)
\end{bmatrix}.
\end{equation}
In addition, the rates of change of the Euler parameters are related to the angular velocity $\om{A}{B}$ as
\begin{equation}\label{eq:QuatRate}
\begin{array}{lclclcl}
\dot{\e{}}_1 &=& \phantom{-} \frac{1}{2} \left( \phantom{-} \eta \w{1} - \e{3} \w{2} + \e{2} \w{3} \right) & , & 
\dot{\e{}}_2 &=& \phantom{-} \frac{1}{2} \left( \phantom{-} \e{3} \w{1} + \eta \w{2} - \e{1} \w{3} \right), \\
\dot{\e{}}_3 &=& \phantom{-} \frac{1}{2} \left(- \e{2} \w{1} + \e{1} \w{2} + \eta \w{3} \right) & , & 
\dot{\eta}    &=& -\frac{1}{2} \left( \phantom{-} \e{1}\w{1} + \e{2}\w{2} + \e{3}\w{3} \right),
\end{array}
\end{equation}
where $\br{\C{B}}{\om{A}{B}} = [\w{1}~\w{2}~\w{3}]\tr$.  Equivalently, the angular velocity components are expressed in terms of the Euler parameters and their rates of change as
\begin{equation}\label{eq:QuatRate-2-Om}
\begin{array}{lcl}
\w{1} &=& 2 \left( \eta \dot{\e{}}_1 - \dot{\eta}\e{1} + \e{3}\dot{\e{}}_2 - \dot{\e{}}_3\e{2} \right), \\
\w{2} &=& 2 \left( \eta \dot{\e{}}_2 - \dot{\eta}\e{2} - \e{3}\dot{\e{}}_1 + \dot{\e{}}_3\e{1} \right), \\
\w{3} &=& 2 \left( \eta \dot{\e{}}_3 - \dot{\eta}\e{3} + \e{2}\dot{\e{}}_1 - \dot{\e{}}_2\e{1} \right).
\end{array}
\end{equation}
The relationships in Eqs.~\eqref{eq:C_BA-Quat}--\eqref{eq:QuatRate-2-Om} will be employed frequently in the derivation of Section~\ref{sect:Derivation1}.

\section{The $rv$-Euler Parameterization \label{sect:Derivation1}}
Consider a particle of mass $m$ located at a point $P$ which moves along with the particle.  Let $\m{r}$ denote the position of $P$ measured relative to an inertially fixed point $O$.  Suppose $O$ is also fixed in frames $\C{E}$ and $\C{A}$ and point $P$ is fixed in frame $\C{B}$.  Here, frame $\C{E}$ is the observation frame in which the relative velocity, $^\C{E}{\m{v}}$, is sought.  In addition, the position and velocity frames, $\C{A}$ and $\C{B}$, are chosen such that the basis vectors $\h{a}_1$ and $\h{b}_1$ are aligned with $\m{r}$ and $^\C{E}{\m{v}}$ respectively.  The remaining degree of freedom for frame $\C{A}$ (rotation of $\h{a}_2$ and $\h{a}_3$ about $\m{r}$) and the remaining degree of freedom for frame $\C{B}$ (rotation of $\h{b}_2$ and $\h{b}_3$ about $^{\C{E}}\m{v}$) are removed by choosing any admissible initial orientation for the remaining basis vectors and applying the nonholonomic constraints $\w{A1}$ = $\w{B1}$ = 0, where $\br{A}{\om{E}{A}} = [\w{A1}~\w{A2}~\w{A3}]\tr$ and $\br{B}{\om{A}{B}} = [\w{B1}~\w{B2}~\w{B3}]\tr$.  The precise nature of the angular velocity constraints are discussed later in the derivation.

Given that the initial orientations for the basis vectors of frames $\C{A}$ and $\C{B}$ are somewhat arbitrary, suppose that choices for the initial orientations of $\{ \h{a}_2, \h{a}_3 \}$ about $\m{r}$ and $\{ \h{b}_2, \h{b}_3 \}$ about $^{\C{E}}\m{v}$ have been made.  Suppose further that the basis vectors assigned to the inertial frame $\C{N}$ and the observation frame $\C{E}$ have appropriate definitions per the application.  For example, the basis vectors associated with the ECI (Earth-centered inertial) and ECEF (Earth-centered, Earth-fixed) frames may be a convenient choice for modeling Earth-relative motion.  After defining the initial orientations of the basis vectors assigned to each frame, the initial values for the Euler parameters $\{ \e{A1},\e{A2},\e{A3},\eta_A \}$ and $\{ \e{B1},\e{B2},\e{B3},\eta_B \}$ are chosen such that the direction cosine matrices
\begin{equation}\label{eq:CAE-final}
\CAE = 
\begin{bmatrix} 
1 - 2(\e{A2}^2 + \e{A3}^2) & 
2(\e{A1} \e{A2} + \e{A3} \eta_A) &
2(\e{A1} \e{A3} - \e{A2} \eta_A) \\
2(\e{A2} \e{A1} - \e{A3} \eta_A) &
1 - 2(\e{A3}^2 + \e{A1}^2) &
2(\e{A2} \e{A3} + \e{A1} \eta_A) \\
2(\e{A3} \e{A1} + \e{A2} \eta_A) &
2(\e{A3} \e{A2} - \e{A1} \eta_A) &
1 - 2(\e{A1}^2 + \e{A2}^2)
\end{bmatrix}, \vspace{2mm}\\
\end{equation}
and
\begin{equation}\label{eq:CBA-final}
\CBA = 
\begin{bmatrix} 
1 - 2(\e{B2}^2 + \e{B3}^2) & 
2(\e{B1} \e{B2} + \e{B3} \eta_B) &
2(\e{B1} \e{B3} - \e{B2} \eta_B) \\
2(\e{B2} \e{B1} - \e{B3} \eta_B) &
1 - 2(\e{B3}^2 + \e{B1}^2) &
2(\e{B2} \e{B3} + \e{B1} \eta_B) \\
2(\e{B3} \e{B1} + \e{B2} \eta_B) &
2(\e{B3} \e{B2} - \e{B1} \eta_B) &
1 - 2(\e{B1}^2 + \e{B2}^2)
\end{bmatrix},
\end{equation}
match the orientations of the relevant basis vectors.  The aforementioned Euler parameters along with the magnitudes $r$ and $v$ of vectors $\m{r}$ and $^\C{E}{\m{v}}$ comprise a ten parameter set termed here as $rv$-Euler parameters.  The differential equations of motion describing the evolution of the $rv$-Euler parameters are derived next.

\subsection{Kinematic Equations\label{sect:Kinematics1}}
The position and velocity of point $P$ relative to point $O$ as viewed by an observer in frame $\C{E}$ is described by
\begin{equation}\label{eq:rv}
\begin{array}{rcl}
\m{r} &=& r \h{a}_1, \\
^\C{E}\m{v} &=& v \h{b}_1,
\end{array}
\end{equation}
where $r$ is the magnitude of the position and $v$ is the speed of the particle as viewed by an observer in the observation frame $\C{E}$.  The relative velocity, $^\C{E}\m{v}$, is also defined by
\begin{equation}\label{eq:kinematics-1}
^\C{E}\m{v} = {\ddtm{\C{E}}{\m{r}}} = {\ddtm{\C{A}}{\m{r}}} + {\om{E}{A}} \times \m{r},
\end{equation}
where $\br{\C{A}}{\om{E}{A}} = [\w{A1}~\w{A2}~\w{A3}]\tr$.  Expressing Eq.~\eqref{eq:kinematics-1} in the basis of frame $\C{A}$ yields
\begin{equation}\label{eq:kinematics-2}
\CAB \br{\C{B}}{^\C{E}\m{v}} = \br{\C{A}}{\ddtm{\C{A}}{\m{r}}}  + \brsk{\C{A}}{\om{E}{A}}  \{ \m{r} \}_{\C{A}},
\end{equation}
which is rewritten in matrix form and simplified to produce
\begin{equation}\label{eq:kinematics-final-0}
\begin{bmatrix}
\dot{r} \\ \phantom{-} r \w{A3} \\ -r \w{A2}
\end{bmatrix}
=
\CAB
\begin{bmatrix}
v \\ 0 \\ 0
\end{bmatrix}.
\end{equation}
Applying the definition of $\CAB$ and solving for $\dot{r}$, $\w{A2}$, and $\w{A3}$ yields the system of equations
\begin{equation}\label{eq:kinematics-final-1}
\begin{array}{rcl}
\dot{r}                   &=& v \left( 1 - 2(\e{B2}^2 + \e{B3}^2) \right), \\
\w{A2} &=& \frac{2 v}{r} \left( \eta_B \e{B2} - \e{B1} \e{B3} \right), \\
\w{A3} &=& \frac{2 v}{r} \left( \eta_B \e{B3} + \e{B1} \e{B2} \right).
\end{array}
\end{equation}
It is apparent from Eqs.~\eqref{eq:kinematics-final-0} and \eqref{eq:kinematics-final-1} that $\w{A1}$ is unconstrained.  In fact, the freedom to choose $\w{A1}$ is a direct consequence of the remaining degree of freedom in defining frame $\C{A}$ (rotations of $\{ \h{a}_2, \h{a}_3 \}$ about $\m{r}$) discussed earlier.  Thus, a whole family of parameterizations exists, each with distinct properties, by simply constraining $\w{A1}$ in various ways.  While it is possible to constrain $\w{A1}$ in an infinite number of ways, the $rv$-Euler parameterization employs $\w{A1} = 0$.  The $\w{A1} = 0$ constraint is equivalently expressed by the Euler parameters and their rates of change as
\begin{equation}\label{eq:w_A1=0}
0 = \eta_A \edot{A1} - \dot{\eta}_A \e{A1} + \e{A3} \edot{A2} - \edot{A3} \e{A2},
\end{equation}
where Eq.~\eqref{eq:QuatRate-2-Om} has been applied.  Thus, by employing $\w{A1} = 0$, the Euler parameter rates are determined by Eq.~\eqref{eq:QuatRate} as
\begin{equation}\label{eq:kinematics-final-2}
\begin{array}{lclclcl}
\edot{A1} &=& -\frac{1}{2} \w{A2} \e{A3} +\frac{1}{2} \w{A3} \e{A2} & , & 
\edot{A2} &=& \phantom{-} \frac{1}{2} \w{A2} \eta_A - \frac{1}{2} \w{A3} \e{A1},  \\
\edot{A3} &=& \phantom{-} \frac{1}{2} \w{A2} \e{A1}  + \frac{1}{2} \w{A3} \eta_A & , & 
\dot{\eta}_A &=& -\frac{1}{2} \w{A2} \e{A2} -\frac{1}{2} \w{A3} \e{A3}.
\end{array}
\end{equation}
Together, Eqs.~\eqref{eq:kinematics-final-1} and \eqref{eq:kinematics-final-2} define the kinematic equations.  Although there are a total of five differential equations which define the evolution of the parameters $\{ r, \e{A1}, \e{A2}, \e{A3}, \eta_{A} \}$, there are only three degrees of freedom because the constraints $\e{A1}^2 + \e{A2}^2 + \e{A3}^2 + \eta_A^2 = 1$ and $\w{A1} = 0$ are both implicit in the equations.

\subsection{Kinetic Equations\label{sect:Dynamics1}}
Now consider the particle kinetics.  Newton's second law is expanded in terms of the relative velocity as
\begin{equation}\label{eq:dynamics-1}
\frac{1}{m}\m{F} = {\ddtm{\C{B}}{^\C{E}\m{v}}} 
+ \left( \om{A}{B} + \om{E}{A} + 2 {\om{N}{E}} \right) \times {^\C{E}\m{v}}
+ \al{N}{E} \times \m{r}
+ \om{N}{E} \times \om{N}{E} \times \m{r},
\end{equation}
where $\br{\C{B}}{\om{A}{B}} = [\w{B1}~\w{B2}~\w{B3}]\tr$, $\br{\C{A}}{\om{E}{A}} = [\w{A1}~\w{A2}~\w{A3}]\tr$, $\br{\C{E}}{\om{N}{E}} = [\w{E1}~\w{E2}~\w{E3}]\tr$, and $\br{\C{E}}{\al{N}{E}} = [\alpha_{E1}~\alpha_{E2}~\alpha_{E3}]\tr$.  Next, let the apparent force, denoted $\tilde{\m{F}}$, be defined as
\begin{equation}\label{eq:F-apparent-1}
\tilde{\m{F}} = \m{F} 
- m \left(
2 {\om{N}{E}} \times {^\C{E}\m{v}}
+ \al{N}{E} \times \m{r}
+ \om{N}{E} \times \om{N}{E} \times \m{r}
\right),
\end{equation}
such that $\frac{1}{m} \tilde{\m{F}} = {^\C{E}{\m{a}}}$.  Coordinatizing Eq.~\eqref{eq:F-apparent-1} in the basis of frame $\C{B}$ and dividing by $m$ yields
\begin{equation}\label{eq:F-apparent-2}
\begin{array}{rcl}
\frac{1}{m}\br{\C{B}}{\tilde{\m{F}}} &\hspace{-2mm} =& \hspace{-2mm} \frac{1}{m} \br{\C{B}}{\m{F}}
- 2 \CBE \brsk{\C{E}}{\om{N}{E}} \CEB \br{\C{B}}{^\C{E}\m{v}}
- \CBE \brsk{\C{E}}{\al{N}{E}} \CEA \br{\C{A}}{\m{r}} \\
& & 
- \CBE \brsk{\C{E}}{\om{N}{E}} \brsk{\C{E}}{\om{N}{E}} \CEA \br{\C{A}}{\m{r}},
\end{array}
\end{equation}
which is expressed in matrix form as
\begin{equation}\label{eq:F-apparent-3}
\begin{array}{rcl}
\frac{1}{m} \begin{bmatrix}
\tilde{f}_1 \\ \tilde{f}_2 \\ \tilde{f}_3
\end{bmatrix}
&=&
\frac{1}{m} \begin{bmatrix}
f_1\\ f_2 \\ f_3
\end{bmatrix}

- 2 v
\begin{bmatrix}
\m{0}\tr \\
\phantom{-}\CBE(3,:) \\
-\CBE(2,:)
\end{bmatrix}
\begin{bmatrix}
\w{E1} \\ \w{E2} \\ \w{E3}
\end{bmatrix}

- r \CBA
\begin{bmatrix}
\m{0}\tr \\
\phantom{-}\CAE(3,:) \\
-\CAE(2,:)
\end{bmatrix}
\begin{bmatrix}
\alpha_{E1} \\ \alpha_{E2} \\ \alpha_{E3}
\end{bmatrix}

\vspace{2mm} \\ & &

- r \CBE
\begin{bmatrix}
-\w{E2}^2 - \w{E3}^2  & \w{E1}\w{E2}              & \w{E1}\w{E3} \\
\w{E2}\w{E1}               & -\w{E1}^2 - \w{E3}^2 & \w{E2}\w{E3} \\
\w{E3}\w{E1}               & \w{E3}\w{E2}              & -\w{E1}^2 - \w{E2}^2
\end{bmatrix}
\CEA(:,1),
\end{array}
\end{equation}
where $\br{\C{B}}{\tilde{\m{F}}} = [\tilde{f}_1~\tilde{f}_2~\tilde{f}_3]\tr$ and $\br{\C{B}}{\m{F}} = [f_1~f_2~f_3]\tr$.  Individual expressions for $\tilde{f}_1$, $\tilde{f}_2$, and $\tilde{f}_3$ in terms of the Euler parameters could be obtained from Eq.~\eqref{eq:F-apparent-3} by appropriate application of $\CAE = \CEA\tr$ in Eq.~\eqref{eq:CAE-final}, $\CBA = \CAB\tr$ in Eq.~\eqref{eq:CBA-final}, and $\CBE = \CBA \CAE$.  However, such an expression is quite lengthy and is omitted for clarity and brevity.

Now return to Eq.~\eqref{eq:dynamics-1}.  Substituting the definition for $\tilde{\m{F}}$, expressing all vectors in the basis of frame $\C{B}$, and reordering terms produces
\begin{equation}\label{eq:dynamics-2}
\br{\C{B}}{\ddtm{\C{B}}{^\C{E}\m{v}}} 
+ \brsk{\C{B}}{\om{A}{B}} \br{\C{B}}{^\C{E}\m{v}}
=
\frac{1}{m} \br{\C{B}}{\tilde{\m{F}}}
- \CBA \brsk{\C{A}}{\om{E}{A}} \CAB \br{\C{B}}{^\C{E}\m{v}},
\end{equation}
which is equivalent to
\begin{equation}\label{eq:dynamics-final-0}
\begin{array}{rcl}
\begin{bmatrix}
\dot{v} \\ \phantom{-} v \w{B3} \\ -v \w{B2}
\end{bmatrix}
&=&
\frac{1}{m} \begin{bmatrix}
\tilde{f}_1 \\ \tilde{f}_2 \\ \tilde{f}_3
\end{bmatrix}

- v
\begin{bmatrix}
\m{0}\tr \\
\phantom{-}\CBA(3,:) \\
-\CBA(2,:)
\end{bmatrix}
\begin{bmatrix}
\w{A1} \\ \w{A2} \\ \w{A3}
\end{bmatrix}.
\end{array}
\end{equation}
Recalling that $\w{A1} = 0$, Eq.~\eqref{eq:dynamics-final-0} is solved for $\dot{v}$, $\w{B2}$, and $\w{B3}$ to produce
\begin{equation}\label{eq:dynamics-final-1}
\begin{array}{rcl}
\dot{v} &=& \frac{1}{m} \tilde{f}_1, \\

\w{B2} &=& - \frac{1}{m v} \tilde{f}_3 
- \w{A2} \left( 1 - 2(\e{B1}^2 + \e{B3}^2) \right)
- 2 \w{A3} (\e{B2} \e{B3} + \e{B1} \eta_B), \\

\w{B3} &=& \frac{1}{m v} \tilde{f}_2 
- 2 \w{A2} (\e{B2} \e{B3} - \e{B1} \eta_B)
- \w{A3} \left( 1 - 2(\e{B1}^2 + \e{B2}^2) \right).
\end{array}
\end{equation}
Similar to the kinematics derivation where $\w{A1}$ was unconstrained, it is observed in Eqs.~\eqref{eq:dynamics-final-0} and \eqref{eq:dynamics-final-1} that $\w{B1}$ is unconstrained.  The freedom to choose $\w{B1}$ stems from the remaining degree of freedom in the definition of frame $\C{B}$ (That is, rotations of $\{ \h{b}_2, \h{b}_3 \}$ about $^{\C{E}}\m{v}$ are arbitrary).  Following the same reasoning as in Section~\ref{sect:Kinematics1}, the constraint $\w{B1} = 0$ is applied to remove the ambiguity in the definition of frame $\C{B}$.  The constraint is equivalent to
\begin{equation}\label{eq:w_B1=0}
0 = \eta_B \edot{B1} - \dot{\eta}_B \e{B1} + \e{B3} \edot{B2} - \edot{B3} \e{B2},
\end{equation}
by Eq.~\eqref{eq:QuatRate-2-Om}.  Thus, the Euler parameter rates are determined by Eq.~\eqref{eq:QuatRate} as
\begin{equation}\label{eq:dynamics-final-2}
\begin{array}{lclclcl}
\edot{B1} &=& -\frac{1}{2} \w{B2} \e{B3} + \frac{1}{2} \w{B3} \e{B2} & , & \edot{B2} &=& \phantom{-} \frac{1}{2} \w{B2} \eta_B - \frac{1}{2} \w{B3} \e{B1}, \\
\edot{B3} &=& \phantom{-} \frac{1}{2} \w{B2} \e{B1} + \frac{1}{2} \w{B3} \eta_B & , & \dot{\eta}_B &=& -\frac{1}{2} \w{B2} \e{B2} - \frac{1}{2} \w{B3} \e{B3}.
\end{array}
\end{equation}
Together, Eqs.~\eqref{eq:F-apparent-3}, \eqref{eq:dynamics-final-1}, and \eqref{eq:dynamics-final-2} define the kinetic equations.  The kinetic equations define the evolution of five parameters, $\{ v, \e{B1}, \e{B2}, \e{B3}, \eta_{B} \}$, in addition to the five introduced in the kinematics derivation, $\{ r, \e{A1}, \e{A2}, \e{A3}, \eta_{A} \}$.  Once again, it is noted that there are only three degrees of freedom, because the constraints $\e{B1}^2 + \e{B2}^2 + \e{B3}^2 + \eta_B^2 = 1$ and $\w{B1} = 0$ are both implicit in the equations.

\section{Examples\label{sect:Examples}}
Next, the equations of motion developed in Section~\ref{sect:Derivation1} are demonstrated on two examples and compared against spherical coordinates.  The first example involves explicit trajectory simulation for a satellite traveling in a circular, sun-synchronous orbit about the Earth.  The example highlights the effect of the singularity at the poles in the spherical parameterization and the lack thereof in the $rv$-Euler parameterization.  Next, the second example involves trajectory optimization via direct collocation for an atmospheric entry maneuver.  The example includes a vertical impact terminal condition that coincides with a singularity in the equations of motion for a spherical coordinate parameterization, whereas the $rv$-Euler parameterization remains well-defined in vertical flight.

\subsection{Example 1: Orbit Propagation}\label{sect:Ex1}
Consider the two-body differential equation \cite{Bate1}
\begin{equation}\label{eq:Ex1-EOM}
\lsup{\C{E}}{\m{a}} = -\frac{\mu_{e}}{r^3} \m{r},
\end{equation}
where frame $\C{E}$ is the Earth-centered inertial (ECI) frame and $\mu_{e}$ is the gravitational parameter of the Earth.  Equation~\eqref{eq:Ex1-EOM} is expressed in spherical coordinates as
\begin{equation}\label{eq:Ex1-EOM-Geo}
\begin{array}{lclclclclcl}
\dot{r} &=& v \sin \gamma & , & 
\dot{\phi} &=& \displaystyle \frac{v}{r \cos \theta} \cos \gamma \sin \psi & , & 
\dot{\theta} &=& \displaystyle \frac{v}{r} \cos \gamma \cos \psi,
\vspace{1mm}\\
\dot{v}  &=&\displaystyle  -\frac{\mu_{e}}{r^{2}} \sin\gamma  & , & 
\dot{\gamma}  &=&\displaystyle  \cos \gamma \left(\frac{v}{r} - \frac{\mu_{e}}{r^2 v} \right) & , & 
\dot{\psi} &=& \displaystyle \frac{v}{r \cos \theta} \cos \gamma \sin \psi  \sin \theta,
\end{array}
\end{equation}
where $r$ is the geocentric radius, $\phi$ is the inertial longitude, $\theta$ is the geocentric latitude, $v$ is the inertial speed, $\gamma$ is the inertial flight path angle, and $\psi$ is the inertial azimuth angle.  Likewise, the $rv$-Euler parameterization of Eq.~\eqref{eq:Ex1-EOM} is given by
\begin{equation}\label{eq:Ex1-EOM-rvEuler}
\begin{array}{rclcrcl}
\dot{r} 
&=& 
v \left( 1 - 2(\e{B2}^2 + \e{B3}^2) \right),

&\quad&

\dot{v} 
&=& 
-\frac{\mu_{e}}{r^2} \left(1 - 2(\e{B2}^2 + \e{B3}^2)\right), \\

\edot{A1} 
&=&
-\frac{1}{2} \w{A2} \e{A3} + \frac{1}{2} \w{A3} \e{A2},

&\quad&

\edot{B1} 
&=& 
-\frac{1}{2} \w{B2} \e{B3} + \frac{1}{2} \w{B3} \e{B2}, \\

\edot{A2} 
&=& 
\phantom{-} \frac{1}{2} \w{A2} \eta_A - \frac{1}{2} \w{A3} \e{A1},

&\quad&

\edot{B2} 
&=& 
\phantom{-} \frac{1}{2} \w{B2} \eta_B - \frac{1}{2} \w{B3} \e{B1}, \\

\edot{A3} 
&=& 
\phantom{-} \frac{1}{2} \w{A2} \e{A1} + \frac{1}{2} \w{A3} \eta_A,

&\quad&

\edot{B3} 
&=& 
\phantom{-} \frac{1}{2} \w{B2} \e{B1} + \frac{1}{2} \w{B3} \eta_B, \\

\dot{\eta}_A 
&=& 
-\frac{1}{2} \w{A2} \e{A2} - \frac{1}{2} \w{A3} \e{A3},

&\quad&

\dot{\eta}_B 
&=& 
-\frac{1}{2} \w{B2} \e{B2} - \frac{1}{2} \w{B3} \e{B3},
\end{array}
\end{equation}
where
\begin{equation}\label{eq:Ex1-EOM-rvEuler-1}
\begin{array}{rcl}
\w{A2} &=& \frac{2 v}{r} \left( \eta_B \e{B2} - \e{B1} \e{B3} \right), \\

\w{A3} &=& \frac{2 v}{r} \left( \eta_B \e{B3} + \e{B1} \e{B2} \right), \\

\w{B2} &=& \phantom{-} \frac{2 \mu_{e}}{r^2 v} \left(\e{B3} \e{B1} + \e{B2} \eta_B\right)
- \w{A2} \left( 1 - 2(\e{B1}^2 + \e{B3}^2) \right)
- 2 \w{A3} (\e{B2} \e{B3} + \e{B1} \eta_B), \\

\w{B3} &=& - \frac{2 \mu_{e}}{r^2 v} \left(\e{B2} \e{B1} - \e{B3} \eta_B\right)
- \w{A3} \left( 1 - 2(\e{B1}^2 + \e{B2}^2) \right)
- 2 \w{A2} (\e{B2} \e{B3} - \e{B1} \eta_B).
\end{array}
\end{equation}

Solutions to Eqs.~\eqref{eq:Ex1-EOM-Geo}--\eqref{eq:Ex1-EOM-rvEuler-1} are obtained on the time interval $t \in [0,T]$, where $T$ is the orbital period.  The initial conditions are given in Table~\ref{tab:Ex1-IC} and each solution is computed numerically via explicit integration.  In the interest of drawing fair comparisons, the integration is carried out using a classical fourth order Runge Kutta method with a constant step size equal to $T/N$, where $N$ is the total number of time steps.  Results are obtained for integer values of $N$ on the interval $[10^1, 10^5]$.

\begin{table}[h]
  \centering
  \caption{Initial Conditions for Example 1.\label{tab:Ex1-IC}}
  \renewcommand{\baselinestretch}{1}\normalsize\normalfont
 \begin{tabular}{|c|c||c|c|}\hline
    Spherical & Initial    & $rv$-Euler   & Initial \\
    Parameters & Value    & Parameters  & Value \\\hline
    
    $r$ & $6971~\textrm{km}$ &
    $r$ & $6971~\textrm{km}$ \\
    
    $\phi$    & $0\deg$ &
    $\e{A1}$ & $-0.753$ \\
    
    $\theta$ & $0\deg$ &
    $\e{A2}$ & $0$ \\
    
    $-$         & $-$ &
    $\e{A3}$ & $0$ \\
    
    $-$        &  $-$ &
    $\eta_A$ & $0.658$ \\
    
    $v$ & $7.562~\textrm{km}/\textrm{s}$ &
    $v$ & $7.562~\textrm{km}/\textrm{s}$ \\
    
    $\gamma$ & $0\deg$ &
    $\e{B1}$    & $0$ \\
    
    $\psi$     & $-172.223\deg$ &
    $\e{B2}$ & $0$ \\
    
    $-$         &  $-$ &
    $\e{B3}$ & $0.707$ \\
    
    $-$         &  $-$ &
    $\eta_B$ & $0.707$ \\\hline
  \end{tabular}
\end{table}

\begin{figure}[hbt!]
  \centering
  \begin{tabular}{lr}
  \subfloat[Inertial longitude, $\phi(t)$ vs. time, $t$.]{\includegraphics[width=.475\textwidth]{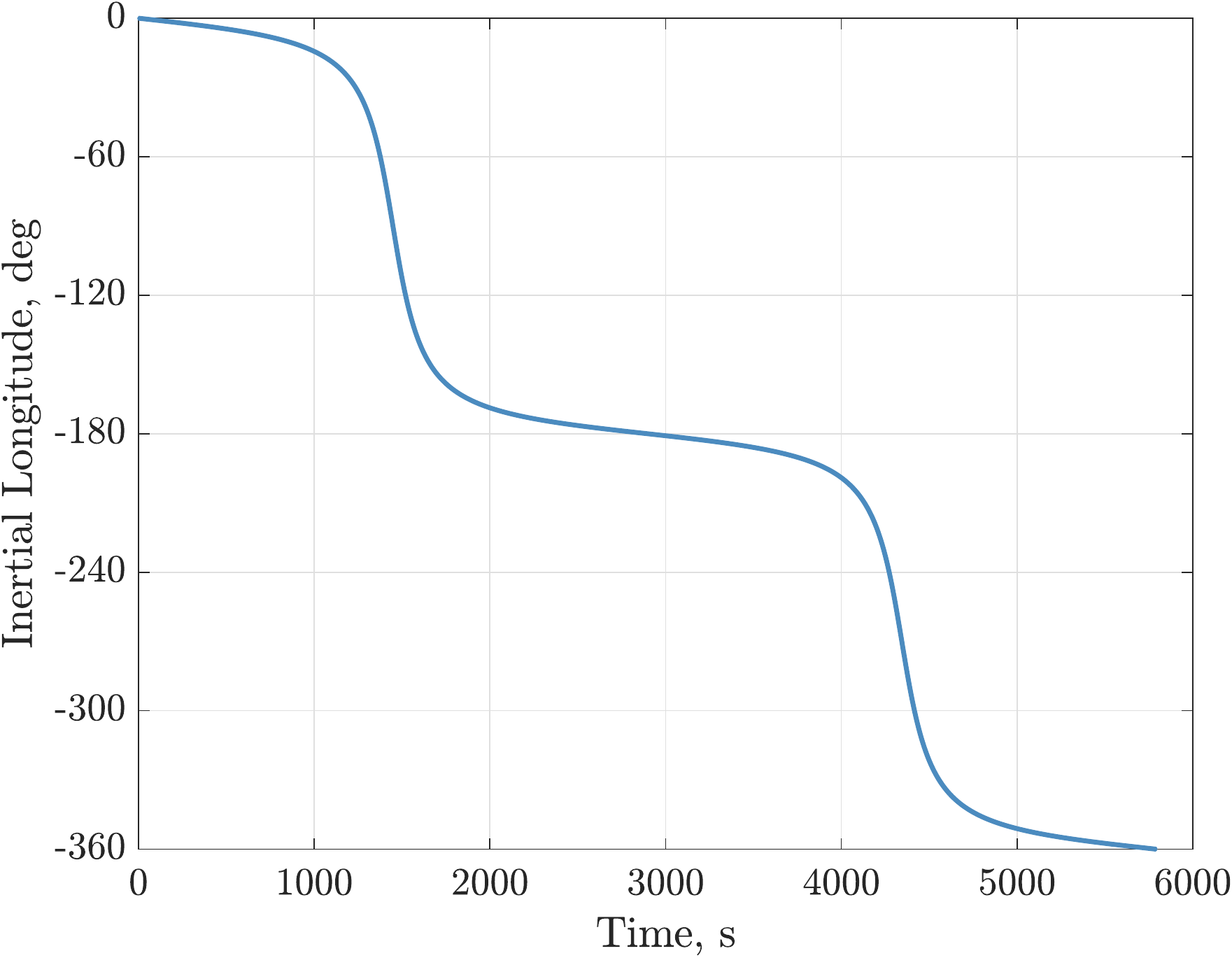}}
  &
  \subfloat[Geocentric latitude, $\theta(t)$ vs. time, $t$.]{\includegraphics[width=.475\textwidth]{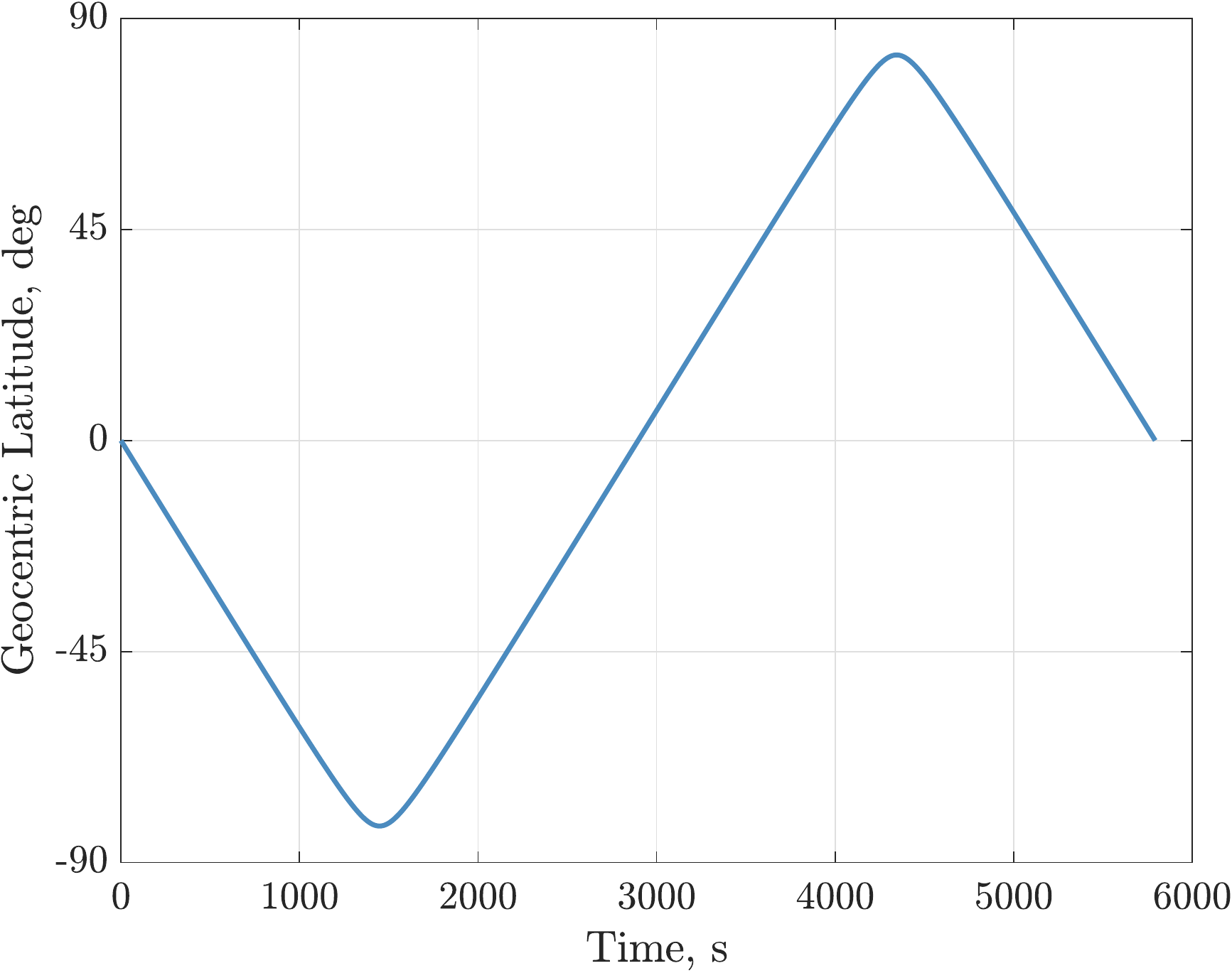}}
  
  \\
  
  \subfloat[Inertial azimuth angle, $\psi(t)$ vs. time, $t$.]{\includegraphics[width=.475\textwidth]{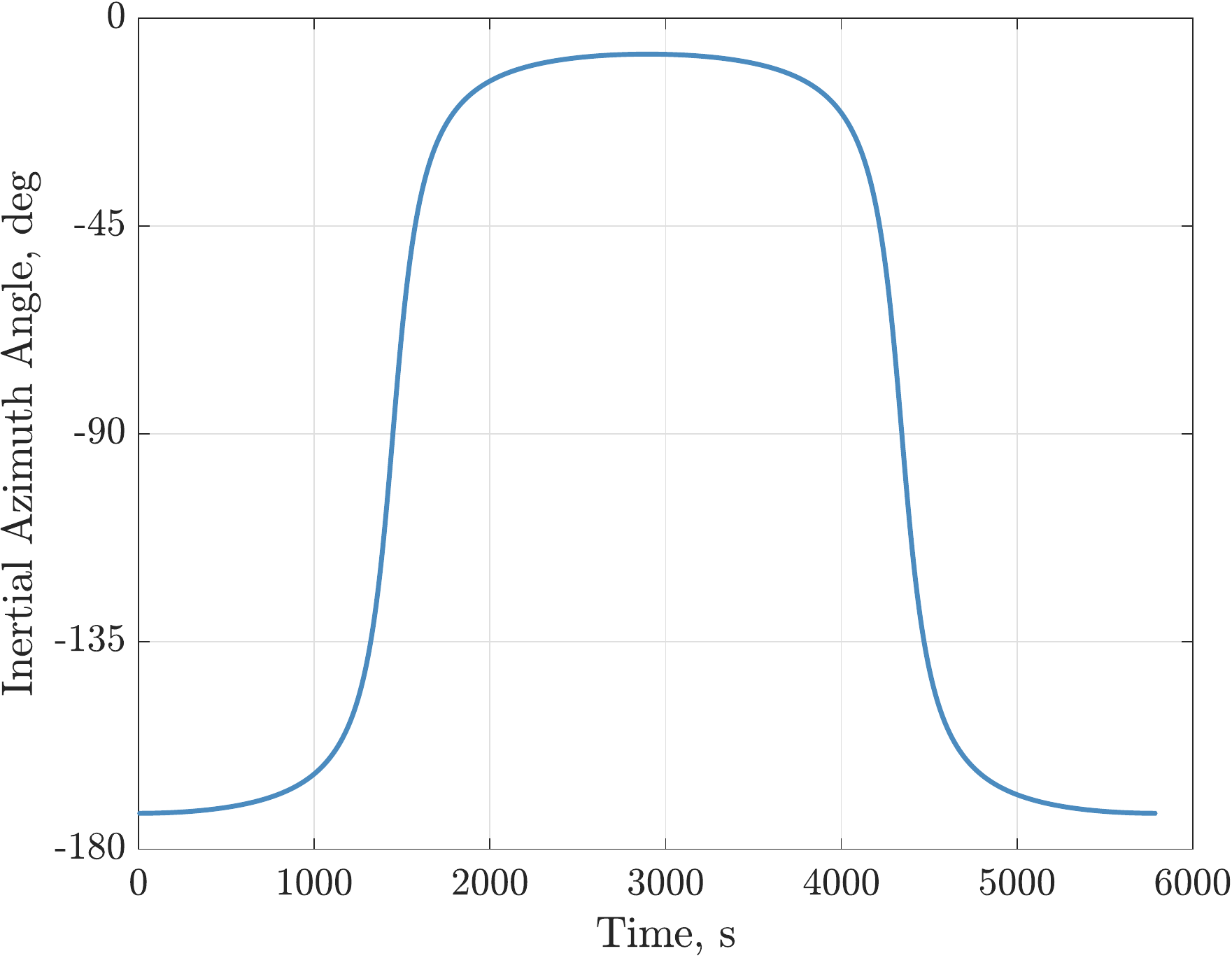}}
  &
  \subfloat[$\CAE$ Euler parameters vs. time, $t$.]{\includegraphics[width=.475\textwidth]{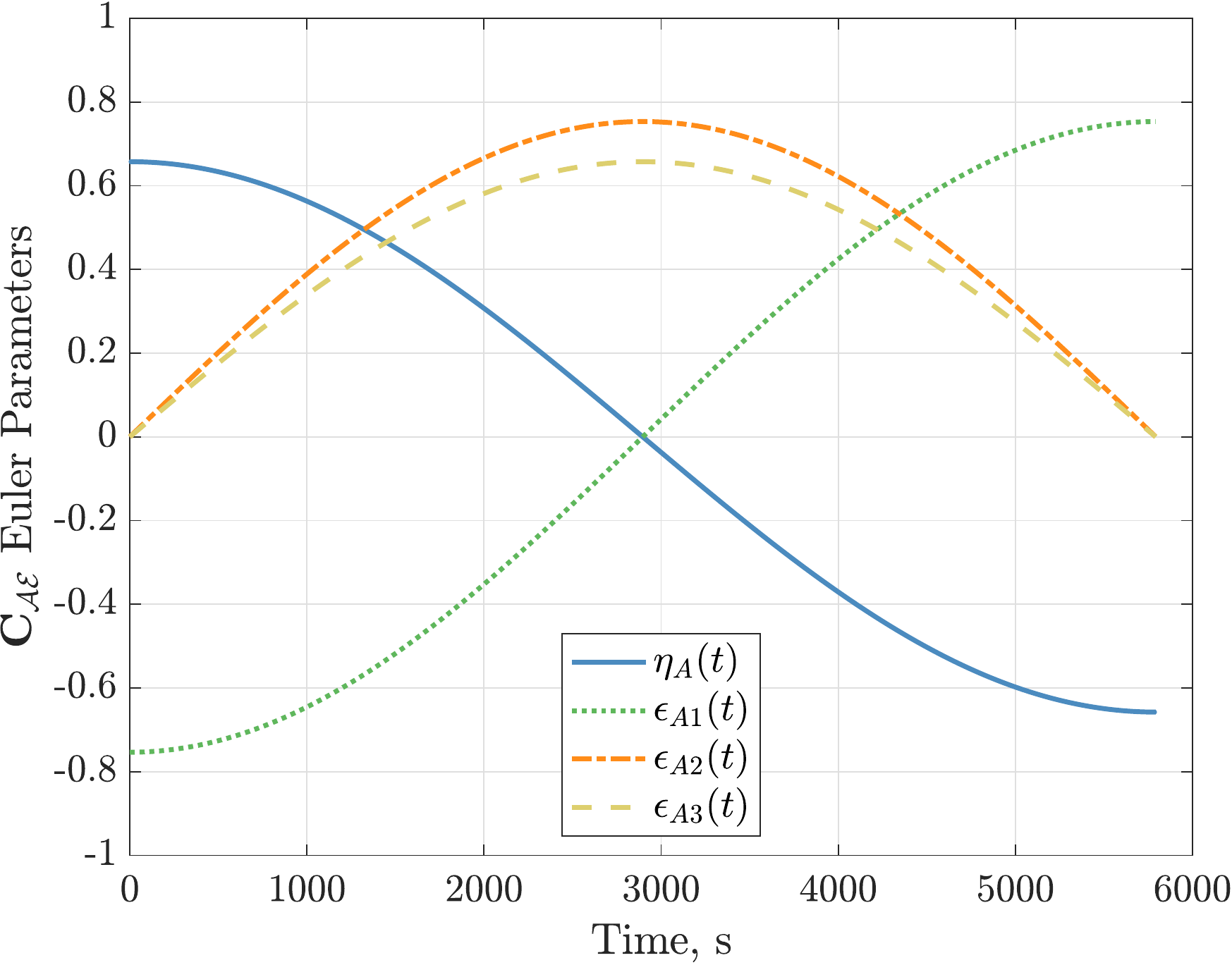}}
  \end{tabular}
  \caption{Simulation results obtained at $N = 10^5$ for Example 1.}\label{fig:Ex1-solution}
\end{figure}

Figure~\ref{fig:Ex1-solution} displays the simulation results obtained at $N = 10^5$ for both the spherical and $rv$-Euler parameterizations.  It is noted that $r$, $v$, $\gamma$, and the $\CBA$ Euler parameters are omitted from Fig.~\ref{fig:Ex1-solution} because they remain essentially constant at their initial values.  Inspection of Fig.~\ref{fig:Ex1-solution} reveals that $\phi$ and $\psi$ change rapidly near the singularities at $\theta \pm 90\deg$.  In contrast, the $\CAE$ Euler parameters smoothly oscillate through one half of a sinusoidal cycle.  The sinusoidal behavior of the $\CAE$ Euler parameters is expected given the half-angle definition of the Euler parameters in Eq.~\eqref{eq:Quat} and noticing that the initial conditions given in Table~\ref{tab:Ex1-IC} align the $\h{a}_3$ and $\h{b}_3$ basis vectors of the position frame $\C{A}$ and the velocity frame $\C{B}$ with the specific angular momentum $^{\C{E}}\m{h} = \m{r} \times {^{\C{E}}\m{v}}$.

After examining Fig.~\ref{fig:Ex1-solution}, it is not immediately clear that the solutions obtained for the spherical and $rv$-Euler parameterizations are in fact identical.  A more quantitative analysis is achieved as follows.  Let the position error, denoted $e_r$, be defined as
\begin{equation}\label{eq:er}
e_r = \left\Vert \br{\C{E}}{\m{r}} - \br{\C{E}}{\m{r}^*} \right\Vert,
\end{equation}
where $\m{r}$ is the position obtained by simulation and $\m{r}^*$ is the analytic solution for the position.  Note that the analytic solution for the position is expressed in the basis of frame $\C{E}$ as
\begin{equation}\label{eq:Ex1-AnalyticTrajectory}
\m{r}^*(t) = (6971~\textrm{km}) 
\left( \cos\left(\frac{2 \pi t}{T}\right) \h{e}_1 + \sin\left(\frac{2 \pi t}{T}\right) (\cos i~\h{e}_2 - \sin i~\h{e}_3) \right), \\
\end{equation}
where $T \approx 5793~\textrm{s}$ is the orbital period and $i \approx 97.8\deg$ is the inclination of the orbit.  

\begin{figure}[hbt!]
  \centering
  \begin{tabular}{lr}
  \subfloat[Position error, $e_r(t)$ vs. time, $t$.]{\includegraphics[width=.475\textwidth]{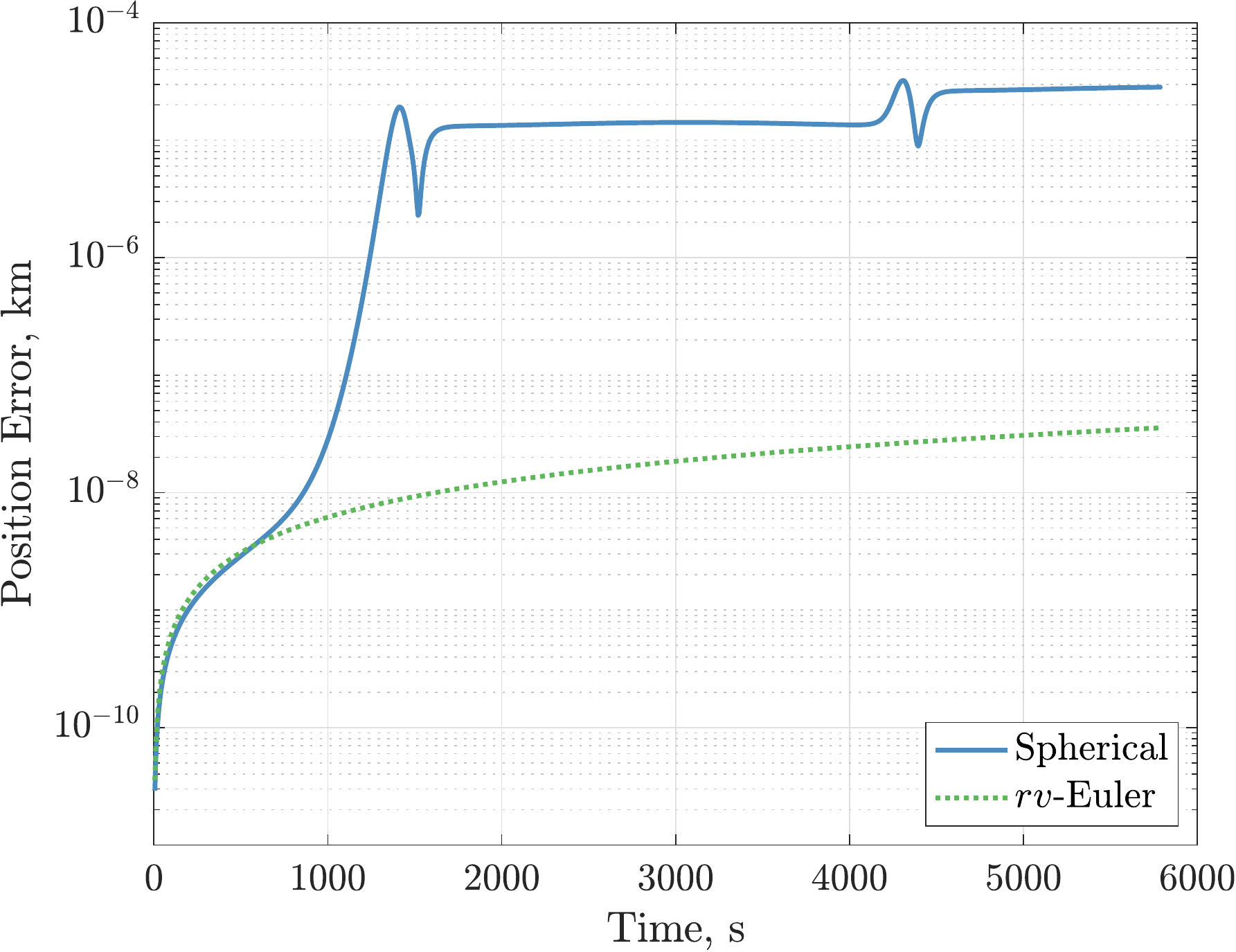}}
  &
  \subfloat[Max position error, $e_{r,\max}$ vs. number of time steps, $N$.]{\includegraphics[width=.475\textwidth]{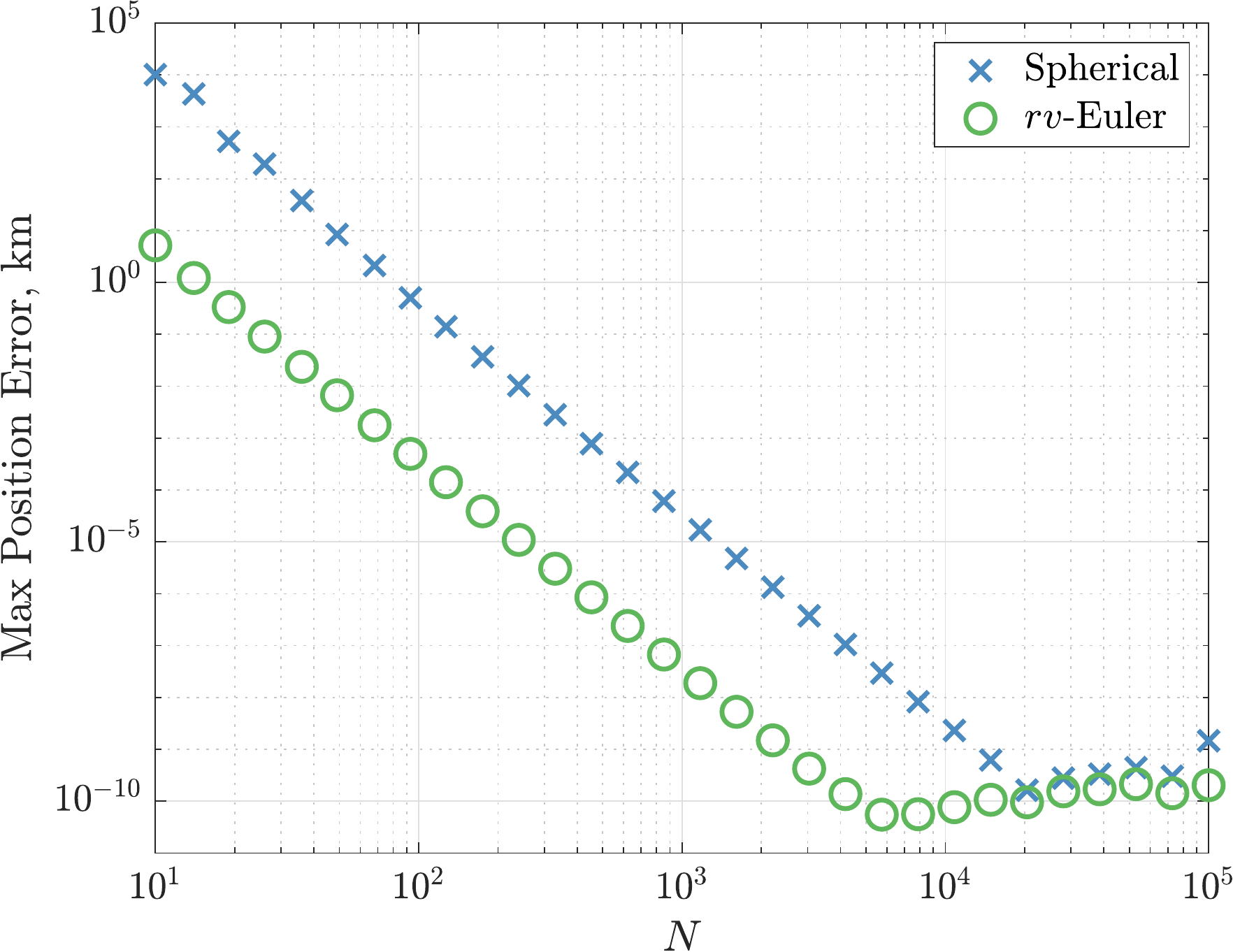}}
  \end{tabular}
  \caption{Position errors at $N = 10^3$ for Example 1 (a).  Maximum position errors for $30$ logarithmically spaced integer values of $N$ spanning $[10^1, 10^5]$ for Example 1 (b).}\label{fig:Ex1-res}
\end{figure}

Applying Eq.~\eqref{eq:er} for each of the simulated trajectories yields the position errors.  Figure~\ref{fig:Ex1-res} (left) illustrates the position errors obtained for both parameterizations when $N = 10^3$ and is representative of the qualitative behavior at most other values of $N$ tested.  Notice how the position error for the spherical parameterization becomes three orders of magnitude larger than that of the $rv$-Euler parameterization after passing near the South Pole singularity.  Figure~\ref{fig:Ex1-res} (right) displays the maximum position error on $t \in [0,T]$, denoted $e_{r,\max}$, for thirty logarithmically spaced integer values of $N$.  It is observed that both parameterizations approach the analytic solution ($e_{r,\max} = 0$) as $N$ increases, but convergence halts around $e_{r,\max} = 10^{-10}$ due to finite computer precision limitations.

\subsection{Example 2: Atmospheric Entry Trajectory Optimization}\label{sect:Ex2}

Consider the following variation of the atmospheric entry optimal control problem detailed in Refs.~\cite{Clarke2,Clarke1}.  The optimal control problems stated in Refs.~\cite{Clarke2,Clarke1} avoid enforcing an exact vertical impact condition (terminal flight path angle equal to $-90\deg$) due to singularity concerns.  Here, the $rv$-Euler parameterization is employed to allay singularity concerns and the vertical impact condition is enforced exactly.

\subsubsection{Problem Statement}\label{sect:Ex2-OCP}

The atmospheric entry optimal problem is stated in terms of the $rv$-Euler parameters as follows.  First, the objective functional to be minimized is given as
\begin{equation}\label{eq:Ex2-Cost}
\C{J} = \int_{0}^{t_f} \left[ 
k_1 \left( \frac{\alpha - \bar{\alpha}}{\alpha_{\max}} \right)^2 + 
k_2 \left( \frac{u_{\alpha}}{u_{\alpha,\max}} \right)^2 +
k_3 \left( \frac{u_{\sigma}}{u_{\sigma,\max}} \right)^2
\right] dt,
\end{equation}
where $\alpha$ is the angle of attack, $u_{\alpha}$ is the angle of attack rate, $u_{\sigma}$ is the bank angle rate, $t_f$ is the terminal time, $(k_1, k_2, k_3)$ are design parameters (and are constants), $\bar{\alpha}$ is the value of $\alpha$ where the lift-drag ratio is a maximum, and $(\alpha_{\max},~u_{\alpha,\max},~u_{\sigma,\max})$ are constants. It is noted that $u_{\alpha}$ and $u_{\sigma}$ are the controls for this example.  Next, the dynamics are given as
\begin{equation}\label{eq:Ex2-EOM}
\begin{array}{rclcrcl}
\dot{r} 
&=& 
v \left( 1 - 2(\e{B2}^2 + \e{B3}^2) \right),

&\quad&

\dot{v} 
&=& 
\frac{1}{m} \tilde{f}_1, \\

\edot{A1} 
&=&
-\frac{1}{2} \w{A2} \e{A3} + \frac{1}{2} \w{A3} \e{A2},

&\quad&

\edot{B1} 
&=& 
-\frac{1}{2} \w{B2} \e{B3} + \frac{1}{2} \w{B3} \e{B2}, \\

\edot{A2} 
&=& 
\phantom{-} \frac{1}{2} \w{A2} \eta_A - \frac{1}{2} \w{A3} \e{A1},

&\quad&

\edot{B2} 
&=& 
\phantom{-} \frac{1}{2} \w{B2} \eta_B - \frac{1}{2} \w{B3} \e{B1}, \\

\edot{A3} 
&=& 
\phantom{-} \frac{1}{2} \w{A2} \e{A1} + \frac{1}{2} \w{A3} \eta_A,

&\quad&

\edot{B3} 
&=& 
\phantom{-} \frac{1}{2} \w{B2} \e{B1} + \frac{1}{2} \w{B3} \eta_B, \\

\dot{\eta}_A 
&=& 
-\frac{1}{2} \w{A2} \e{A2} - \frac{1}{2} \w{A3} \e{A3},

&\quad&

\dot{\eta}_B 
&=& 
-\frac{1}{2} \w{B2} \e{B2} - \frac{1}{2} \w{B3} \e{B3},
\end{array}
\end{equation}
where $m$ is the mass, $r$ is the geocentric radius, $\{ \e{A1}, \e{A2},\e{A3}, \eta_A \}$ are the $\CAE$ Euler parameters, $v$ is the Earth-relative speed, and $\{ \e{B1}, \e{B2},\e{B3}, \eta_B \}$ are the $\CBA$ Euler parameters.  Note that frame $\C{N}$ is the Earth-centered inertial (ECI) frame, frame $\C{E}$ is the Earth-centered, Earth-fixed (ECEF) frame, and frames $\C{A}$ and $\C{B}$ are the position and velocity frames as defined in Section~\ref{sect:Derivation1}.  Moreover, the angular velocity terms in Eq.~\eqref{eq:Ex2-EOM} are given by
\begin{equation}\label{eq:Ex2-EOM-w}
\begin{array}{rcl}
\w{A2} &=& \frac{2 v}{r} \left( \eta_B \e{B2} - \e{B1} \e{B3} \right), \\

\w{A3} &=& \frac{2 v}{r} \left( \eta_B \e{B3} + \e{B1} \e{B2} \right), \\

\w{B2} &=& - \frac{1}{m v} \tilde{f}_3 
- \w{A2} \left( 1 - 2(\e{B1}^2 + \e{B3}^2) \right)
- 2 \w{A3} (\e{B2} \e{B3} + \e{B1} \eta_B), \\

\w{B3} &=& \frac{1}{m v} \tilde{f}_2 
- 2 \w{A2} (\e{B2} \e{B3} - \e{B1} \eta_B)
- \w{A3} \left( 1 - 2(\e{B1}^2 + \e{B2}^2) \right).
\end{array}
\end{equation}
The apparent force terms in Eqs.~\eqref{eq:Ex2-EOM} and \eqref{eq:Ex2-EOM-w} are given as
\begin{equation}\label{eq:Ex2-EOM-Fapp}
\begin{array}{rcl}
\frac{1}{m} \begin{bmatrix}
\tilde{f}_1 \\ \tilde{f}_2 \\ \tilde{f}_3
\end{bmatrix}
&=&
\frac{1}{m} \begin{bmatrix}
f_1 \\ f_2 \\ f_3
\end{bmatrix}

- 2 \w{e} v
\begin{bmatrix}
0 \\ \phantom{-}\CBE(3,3) \\ -\CBE(2,3)
\end{bmatrix}

- r \w{e}^2 \CBA 
\begin{bmatrix}
-\left(\CAE(2,3)\right)^2 - \left(\CAE(3,3)\right)^2 \\
\CAE(2,3) \CAE(1,3) \\
\CAE(3,3) \CAE(1,3)
\end{bmatrix},
\end{array}
\end{equation}
where $\w{e}$ is the rotation rate of the Earth, and $\CAE$ and $\CBA$ are defined in Eqs.~\eqref{eq:CAE-final} and \eqref{eq:CBA-final}.  The forces present are lift, drag, and gravity.  Thus, the net force components in Eq.~\eqref{eq:Ex2-EOM-Fapp} are defined as
\begin{equation}\label{eq:Ex2-EOM-Ftot}
\begin{array}{rcl}
\begin{bmatrix}
{f}_1 \\ {f}_2 \\ {f}_3
\end{bmatrix}
&=&
\begin{bmatrix}
-D \\ L \cos \sigma \\ L \sin \sigma
\end{bmatrix}

-m \frac{\mu_{e}}{r^2} \CBA(:,1),
\end{array}
\end{equation}
where $L=q S C_L(\alpha)$ is the magnitude of the lift force, $D = q S C_D(\alpha)$ is the magnitude of the drag force, $q=\rho v^2/2$ is the dynamic pressure, $S$ is the reference area, $\sigma$ is the bank angle (defined as a rotation about $\h{b}_1$ from the $\h{b}_2$ direction to the lift direction), and $\mu_{e}$ is the gravitational parameter of the Earth.  The details of the aerodynamic model are omitted here but can be found in Refs.~\cite{Clarke2,Clarke1}.

Next, $\alpha$ and $\sigma$ are augmented components of the state such that
\begin{equation}
  \begin{array}{lclclcl}
    \dot{\alpha}& =& u_{\alpha} & \textrm{and} &     \dot{\sigma}& =& u_{\sigma}.
  \end{array}
\end{equation}
Furthermore, the following path constraints are imposed during flight (where $ a_{\max}$, $q_{\min}$, $u_{\alpha,\max}$, and $u_{\sigma,\max}$ are constants):
\begin{equation}\label{eq:Ex2-path}
\begin{array}{rclcrcl}
  \sqrt{L^2+D^2}/m & \leq & a_{\max} & , & q & \geq & q_{\min}, \\
  \alpha & \geq & 0 & , & \alpha & \leq & \alpha_{\max}, \\
  |u_{\alpha}| & \leq & u_{\alpha,\max}  & , & |u_{\sigma}| & \leq & u_{\sigma,\max}.
\end{array}
\end{equation}
Next, the boundary conditions are given in Table~\ref{tab:Ex2-BCs}.  Note that requiring $\e{B1}(t_f) = \eta_B(t_f) = 0$ imposes the vertical impact condition by forcing the velocity direction ($\h{b}_1$) at impact to point in the opposite direction as the position direction ($\h{a}_1$) at impact.  Finally, the endpoint constraints
\begin{equation}\label{eq:Ex2-endpoint}
\begin{array}{lcl}
1 - 2 \left( \e{A2}^2 + \e{A3}^2 \right) \bigr \rvert_{t_f} & = & \cos\theta_T \cos\phi_T, \\
2 \left( \e{A1} \e{A2} + \e{A3} \eta_{A} \right) \bigr \rvert_{t_f} & = & \cos\theta_T \sin\phi_T, \\
2 \left( \e{A1} \e{A3} - \e{A2} \eta_{A} \right) \bigr \rvert_{t_f} & = & \sin\theta_T,
\end{array}
\end{equation}
force the terminal position to coincide with the Earth-relative longitude and geocentric latitude of the target, where $\phi_T = 25.15\deg$ is the Earth-relative longitude of the target and $\theta_T = 0\deg$ is the geocentric latitude of the target.

\begin{table}[h]
  \centering
  \caption{Boundary conditions for Example 2.\label{tab:Ex2-BCs}}
  \renewcommand{\baselinestretch}{1}\normalsize\normalfont
 \begin{tabular}{|c|c||c|c|}\hline
    Symbol & Value & Symbol & Value \\\hline
    $t_0$ & $0~\textrm{s}$ & $t_f$ & FREE \\
    $r(t_0)$         & $R_{e} + 37~\textrm{km}$ & $r(t_f)$ & $R_{e}$ \\
    $\e{A1}(t_0)$ & $0$ & $\e{A1}(t_f)$ & FREE \\
    $\e{A2}(t_0)$ & $0$ & $\e{A2}(t_f)$ & FREE \\
    $\e{A3}(t_0)$ & $0$ & $\e{A3}(t_f)$ & FREE \\
    $\eta_A(t_0)$ & $1$ & $\eta_A(t_f)$ & FREE \\
    $v(t_0)$         & $7.138~\textrm{km/s}$ & $v(t_f)$ & $1.219~\textrm{km/s}$ \\
    $\e{B1}(t_0)$ & $0$ & $\e{B1}(t_f)$ & $0$ \\
    $\e{B2}(t_0)$ & $0$ & $\e{B2}(t_f)$ & FREE \\
    $\e{B3}(t_0)$ & $\sqrt{2} / 2$ & $\e{B3}(t_f)$ & FREE \\
    $\eta_B(t_0)$ & $\sqrt{2} / 2$ & $\eta_B(t_f)$ & $0$ \\
    $\alpha(t_0)$ & FREE & $\alpha(t_f)$ & $0\deg$ \\
    $\sigma(t_0)$ & FREE & $\sigma(t_f)$ & FREE \\\hline
  \end{tabular}
\end{table}

The atmospheric entry optimal control problem is summarized as follows.  Determine the state $( r(t), $ $\e{A1}(t), \e{A2}(t),$ $\e{A3}(t), \eta_A(t), v(t), \e{B1}(t), \e{B2}(t), \e{B3}(t), \eta_B(t), \alpha(t), \sigma(t) )$ and the control $\left( u_{\alpha}(t), u_{\sigma}(t) \right)$ on the time interval $t \in [0,t_f]$ which minimizes the cost functional of Eq.~\eqref{eq:Ex2-Cost} while satisfying the state dynamics, 
the path constraints,
the endpoint constraints,
and the boundary conditions.



\subsubsection{Solution Method and Implementation}\label{sect:Ex2-Implementation}
The atmospheric entry problem described in Section~\ref{sect:Ex2-OCP} is implemented and solved using $hp$-adaptive Gaussian quadrature collocation via the software of Ref.~\cite{Patterson2014}.
Special precautions must be taken when implementing Eq.~\eqref{eq:Ex2-EOM} using a collocation method.  
Recall that the unit norm constraints $\e{A1}^2 + \e{A2}^2 + \e{A3}^2 + \eta_A^2 = 1$ and $\e{B1}^2 + \e{B2}^2 + \e{B3}^2 + \eta_B^2 = 1$, as well as the angular velocity constraints $\w{A1} = 0$ and $\w{B1} = 0$ are all implicit in Eq.~\eqref{eq:Ex2-EOM}.  Thus, the system of differential equations loses four degrees of freedom and may become inconsistent when solving the problem numerically.  The remedy employed here is to introduce four additional control variables, denoted ($u_1$,$u_2$,$u_3$,$u_4$), such that $\edot{A2}  =  u_1$, $ \edot{A3} =  u_2$, $\edot{B2}  =  u_3$, $\edot{B3}  = u_4$
replaces the corresponding differential equations in Eq.~\eqref{eq:Ex2-EOM}.  Next, the additional path constraints
\begin{equation}
\begin{array}{rclcrcl}
| u_1 - \left(
\frac{1}{2} \w{A2} \eta_A - \frac{1}{2} \w{A3} \e{A1} \right) | & \leq & \delta

&,&

| u_2 - \left(
\frac{1}{2} \w{A2} \e{A1} + \frac{1}{2} \w{A3} \eta_A \right) | & \leq & \delta, \\

| u_3 - \left(
\frac{1}{2} \w{B2} \eta_B - \frac{1}{2} \w{B3} \e{B1} \right) | & \leq & \delta

&,&

| u_4 - \left(
\frac{1}{2} \w{B2} \e{B1} + \frac{1}{2} \w{B3} \eta_B \right) | & \leq & \delta,
\end{array}
\end{equation}
are enforced with $\delta = 10^{-6}$ chosen to be one order of magnitude larger than the NLP solver accuracy tolerance.

\subsubsection{Optimal Entry Trajectory}\label{sect:Ex2-solution}
\begin{figure}[hbt!]
  \centering
  \begin{tabular}{lr}
  \subfloat[Altitude, $h(t)$ vs. time, $t$.]{\includegraphics[width=.475\textwidth]{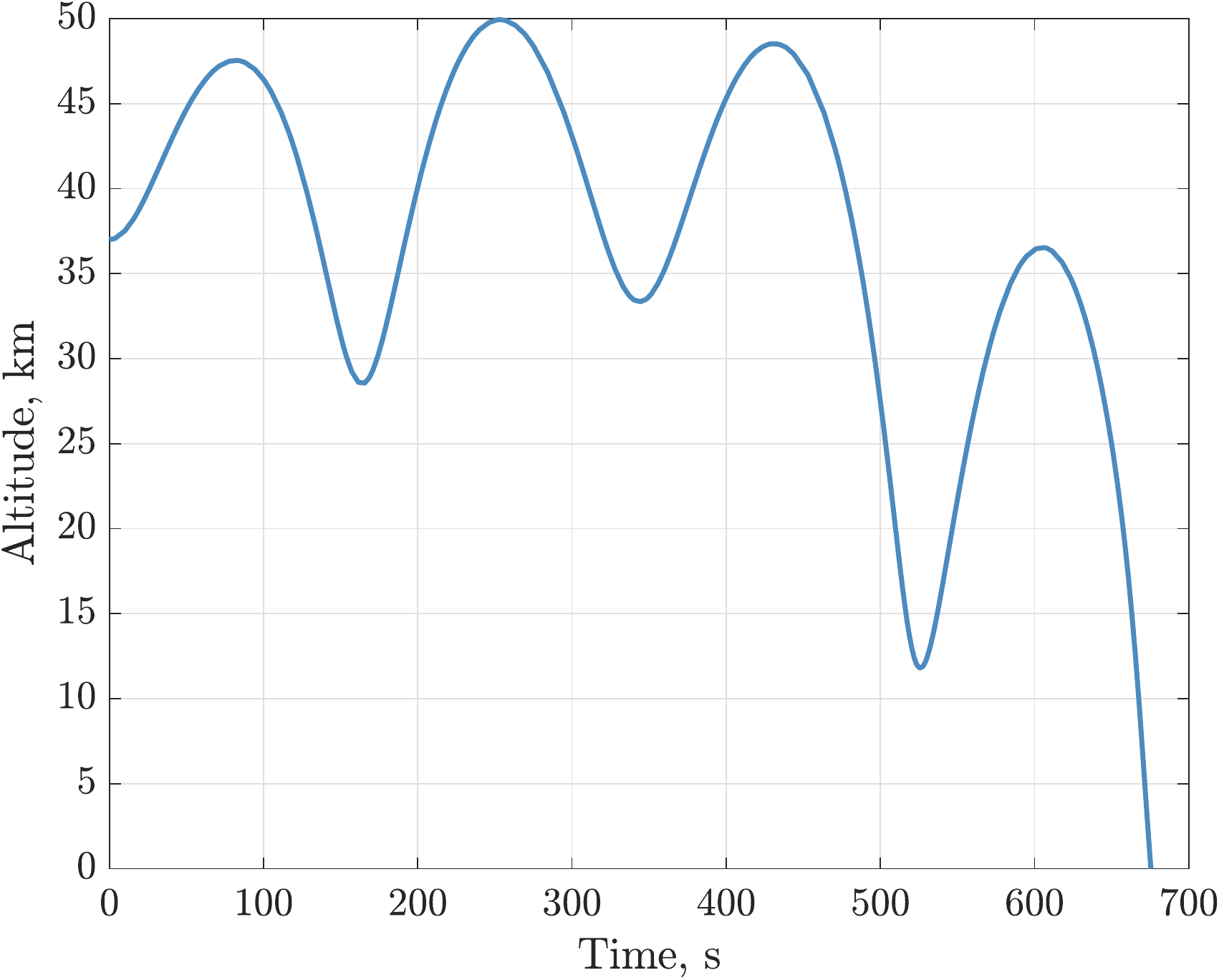}}
  &
  \subfloat[Earth-relative speed, $v(t)$ vs. time, $t$.]{\includegraphics[width=.475\textwidth]{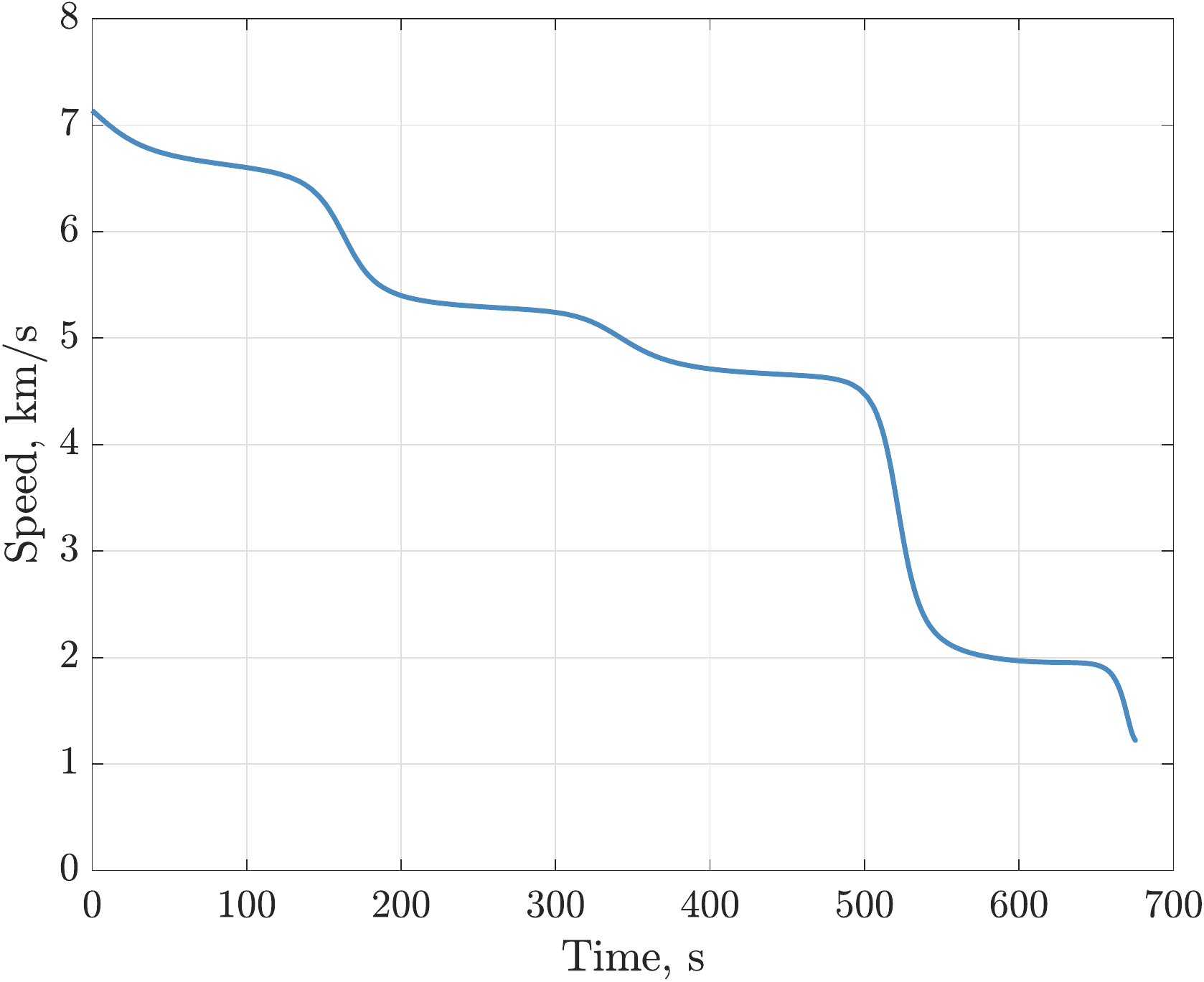}}
  
  \\
  
  \subfloat[$\CAE$ Euler parameters vs. time, $t$.]{\includegraphics[width=.475\textwidth]{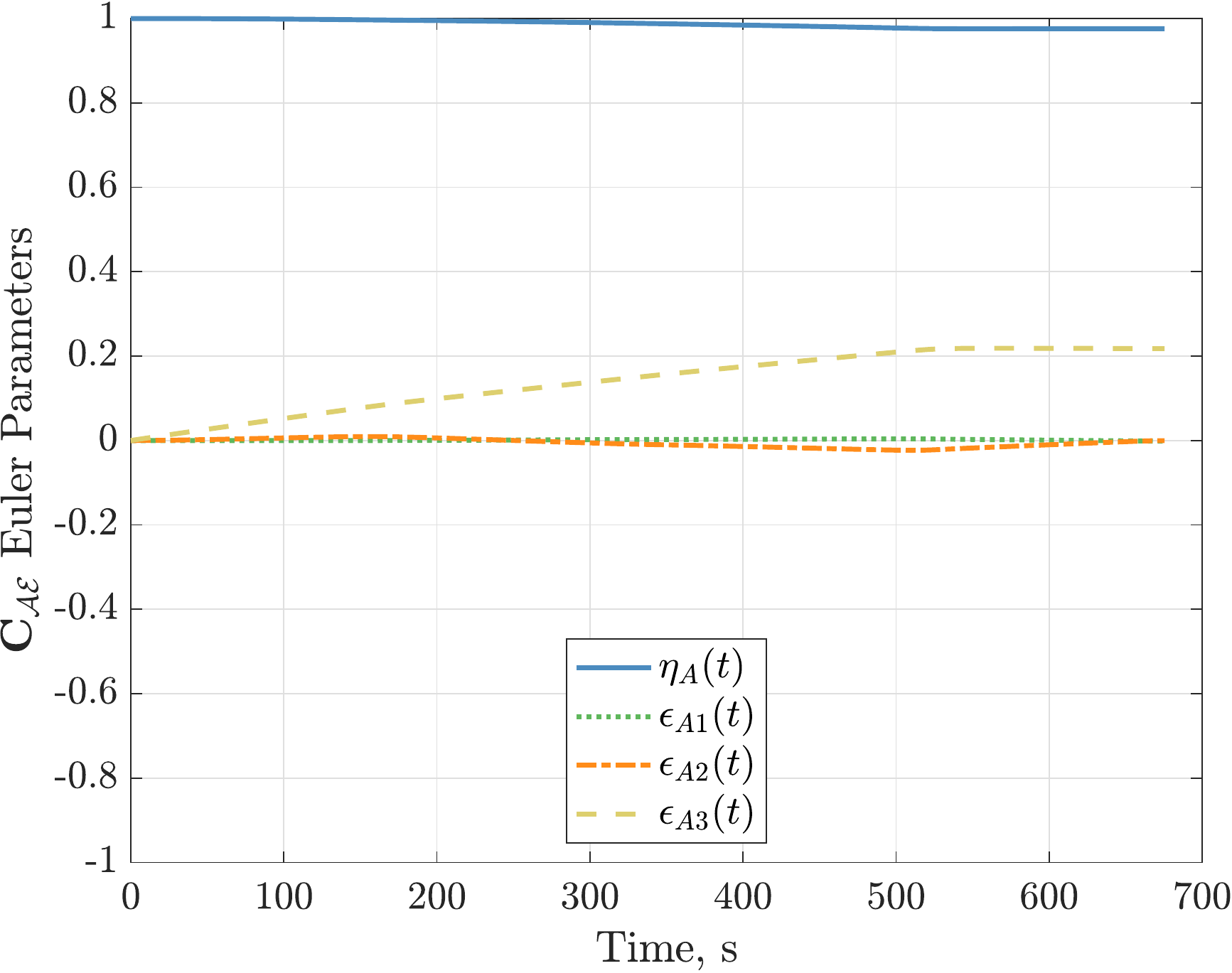}}
  &
  \subfloat[$\CBA$ Euler parameters vs. time, $t$.]{\includegraphics[width=.475\textwidth]{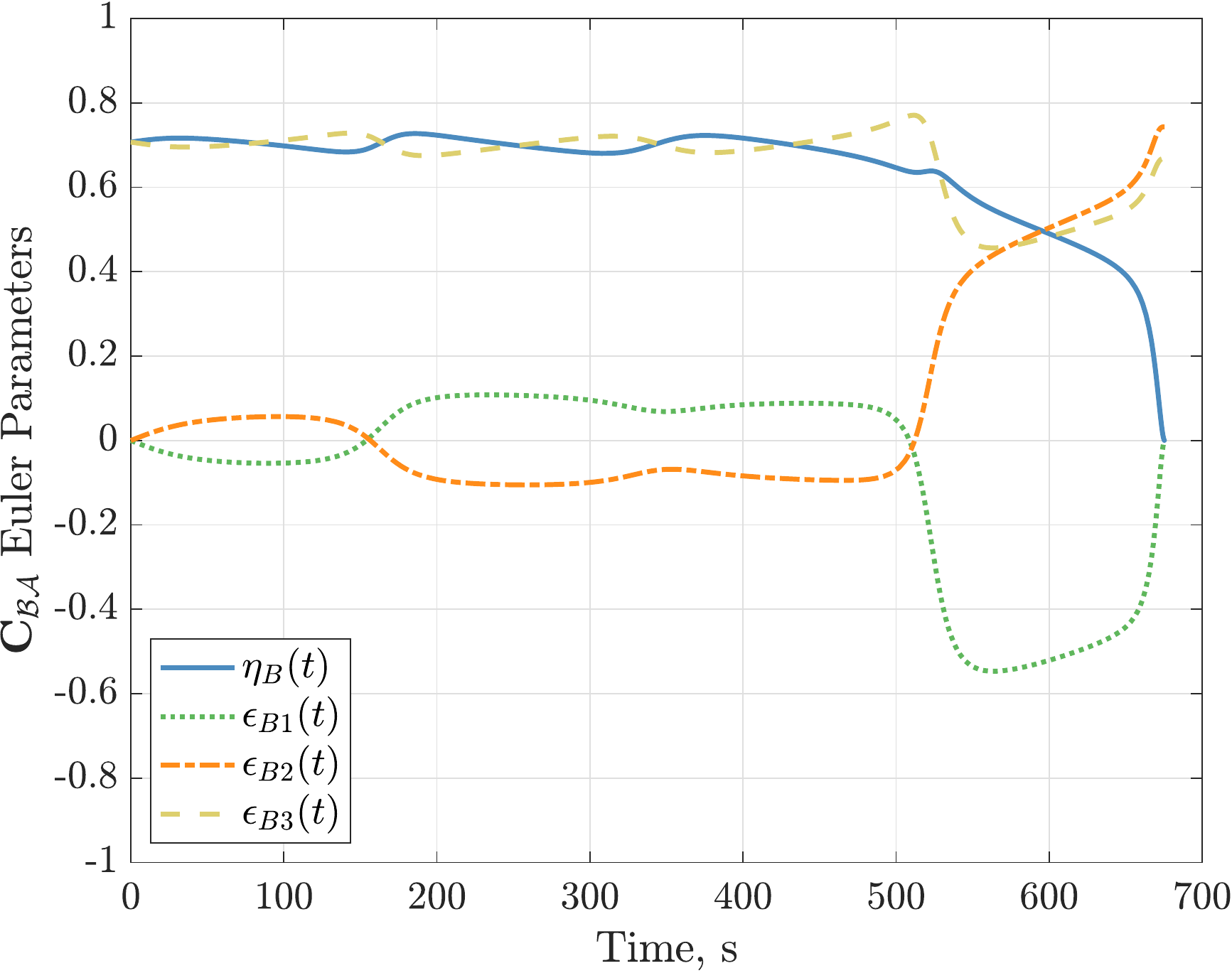}}\label{fig:Ex2-sol-1-CBA}
  
  \\
  
  \subfloat[Angle of attack, $\alpha(t)$ vs. time, $t$.]{\includegraphics[width=.475\textwidth]{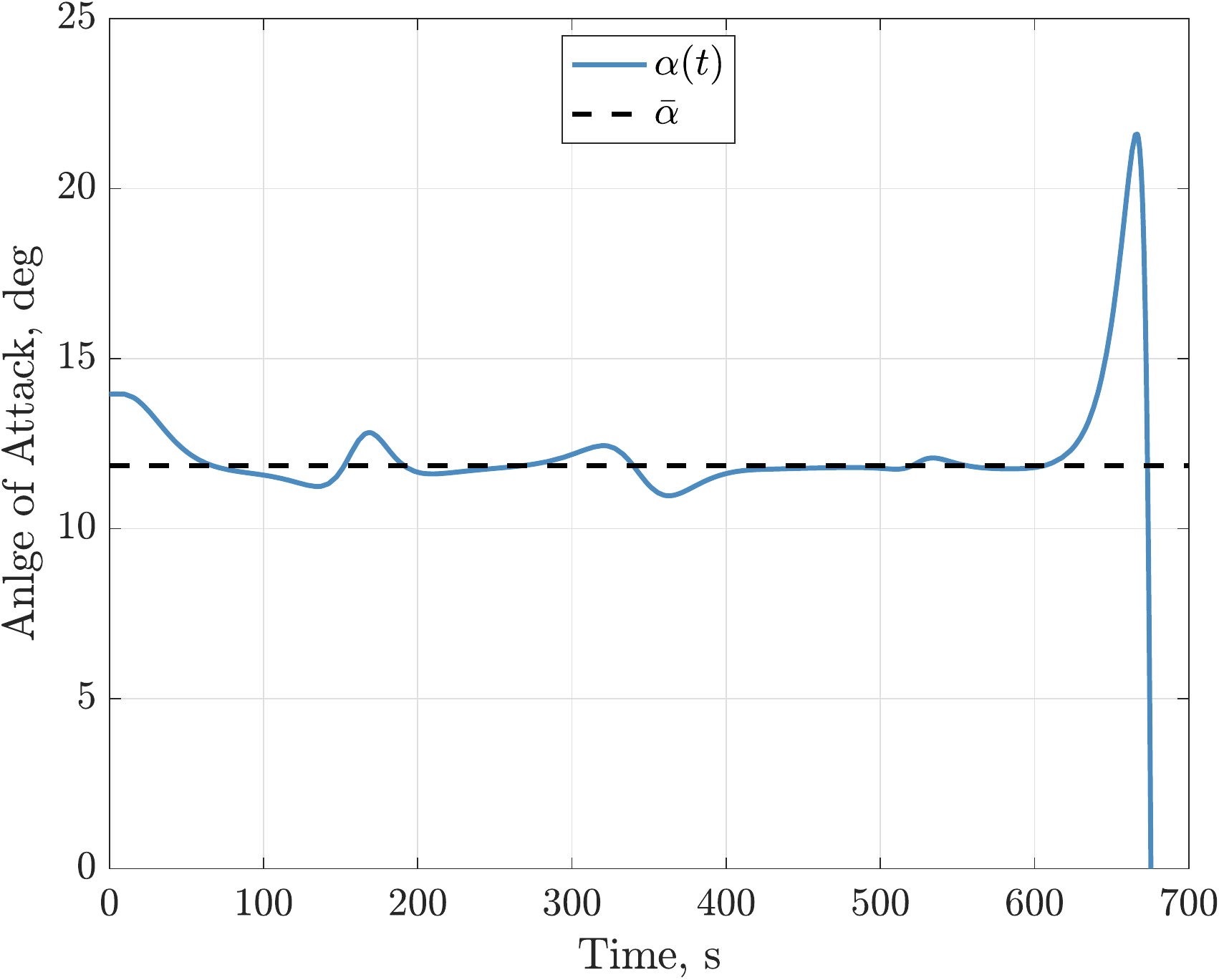}}
  &
  \subfloat[Bank angle, $\sigma(t)$ vs. time, $t$]{\includegraphics[width=.475\textwidth]{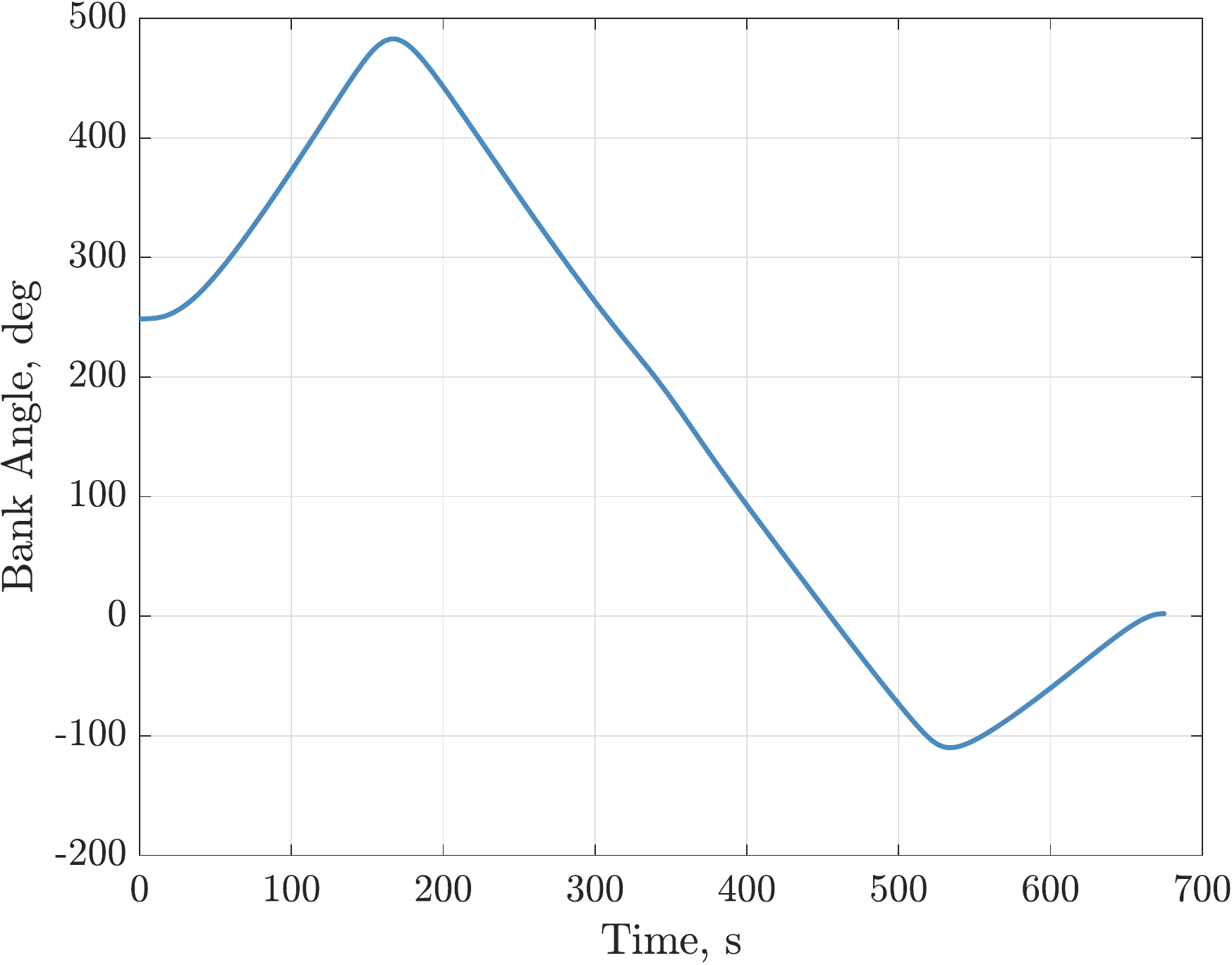}}
  \end{tabular}
  \caption{Optimal solution obtained for Example 2.}\label{fig:Ex2-sol-1}
\end{figure}

Figure~\ref{fig:Ex2-sol-1} shows the optimal entry trajectory where it is seen that the terminal state corresponds to vertically downward flight ($\e{B1}(t_f) = \eta_B(t_f) = 0$).  Notice that Eqs.~\eqref{eq:Ex2-EOM} and \eqref{eq:Ex2-EOM-w} remain well defined at $\e{B1} = \eta_B = 0$.  In contrast, in a spherical coordinate parameterization the differential equation for azimuth, given by
\begin{equation}\label{eq:Ex2-azidot}
\begin{array}{rcl}
\dot{\psi} & =& \frac{L \sin\sigma}{m v \cos\gamma}
+ \frac{v}{r} \cos\gamma \sin\psi \tan\theta
- 2 \w{e} \left( \tan\gamma \cos\psi \cos\theta - \sin\theta \right)
+ \frac{r \w{e}^2}{v \cos\gamma} \sin\psi \sin\theta \cos\theta,
\end{array}
\end{equation}
where $\gamma$ is the Earth-relative flight path angle, is not well-defined for vertical flight ($\gamma = \pm 90$ deg).  Thus, the vertical flight singularity using spherical coordinates is eliminated when using the $rv$-Euler parameterization.

\section{Discussion\label{sect:Discussion}}
The two examples in Section~\ref{sect:Examples} highlight different aspects of the $rv$-Euler parameterization.  The first example validates the accuracy of the $rv$-Euler parameterization by direct comparison of the numerical results with the analytic trajectory.  Moreover, the example contrasts the $rv$-Euler parameters with spherical coordinates and draws attention to the effects of singularities at the North and South poles which are present in the spherical parameterization but not in the $rv$-Euler parameterization.  Specifically, in Example~1 the $rv$-Euler parameterization consistently produces max position errors roughly three orders of magnitude less than the spherical parameterization.  Next, Example 2 studies a more complex problem with no known analytic solution.  This second example focuses on the ability to use the $rv$-Euler parameters to reformulate the trajectory optimization problems in Refs.~\cite{Clarke2,Clarke1} into a form that does not contain a singularity in vertical flight.  

Finally, in addition to the absence of singularities, the $rv$-Euler parameterization does not involve the use of trigonometric functions.  The significance of avoiding the use of trigonometric functions is related to the computational expense in evaluating a trigonometric function.   For example, consider the following software implementation of $\sin x$ \cite{netlib}, given by
 \begin{equation}\label{eq:sin(x)}
 \sin x \approx x + x z (S_1 + z (S_2 + z (S_3 + z (S_4 + z (S_5 +z S_6))))),
 \end{equation}
where $z = x^2$ and $S_i,~i = 1,\ldots,6$ are the coefficients of the truncated power series approximation.  The sine evaluation in Eq.~\eqref{eq:sin(x)} requires $14$ floating-point operations, not to mention the function calls, decision tree, etc. which are also executed in the code.  Now, observe in Eqs.~\eqref{eq:Ex1-EOM-Geo} and \eqref{eq:Ex1-EOM-rvEuler} from Section~\ref{sect:Ex1} that $\sin \gamma$ in the spherical parameterization is equivalently expressed as $(1 - 2(\e{B2}^2 + \e{B3}^2))$ in the $rv$-Euler parameterization.  The reduction in the number of floating-point operations is clear.

\section{Conclusions\label{sect:Conclusion}}
An $rv$-Euler parameterization of the equations of motion for a point mass has been derived.  The equations of motion have been shown to avoid singularities found in commonly used parameterizations.  Moreover, it was found that both position-dependent and velocity-dependent forces could be quantified in a tractable manner using the $rv$-Euler parameters.  Finally, two examples were used to demonstrate the accuracy and nonsingular nature of the equations of motion. 

\section*{Acknowledgments}
The authors gratefully acknowledge support for this research from the from the U.S.~National Science Foundation under grants CMMI-1563225, DMS-1522629, and DMS-1819002, from the U.S.~Office of Naval Research under grant N00014-19-1-2543, and from the U.S.~Department of Defense under the National Defense Science \& Engineering Graduate Fellowship (NDSEG) Program.

\renewcommand{\baselinestretch}{1}
\normalsize\normalfont

\bibliographystyle{aiaa}     

\begin{thebibliography}{10}
\newcommand{\enquote}[1]{``#1''}

\bibitem{Shuster1993}
Shuster, M.~D., \enquote{A Survey of Attitude Representations,} {\em The
  Journal of the Astronautical Sciences\/}, Vol.~41, No.~4, oct 1993,
  pp.~439--517.

\bibitem{Kane1983}
Kane, T.~R., Likins, P.~W., and Levinson, D.~A., {\em Spacecraft Dynamics\/},
  McGraw-Hill Book Co, 1983.

\bibitem{Wie1985}
Wie, B. and Barba, P.~M., \enquote{Quaternion Feedback for Spacecraft Large
  Angle Maneuvers,} {\em Journal of Guidance, Control, and Dynamics\/}, Vol.~8,
  No.~3, may 1985, pp.~360--365. \url{https://doi.org/10.2514/3.19988}.

\bibitem{Wie2008}
Wie, B., {\em Space Vehicle Dynamics and Control\/}, American Institute of
  Aeronautics and Astronautics, 2008.

\bibitem{Hughes2}
Hughes, P.~C., {\em Spacecraft Attitude Dynamics\/}, Dover Publications,
  Mineola, NY, 2012.

\bibitem{Junkins2012}
Junkins, J.~L. and Turner, J.~D., {\em Optimal Spacecraft Rotational
  Maneuvers\/}, Elsevier, 2012.

\bibitem{Wertz2012}
Wertz, J.~R., {\em Spacecraft Attitude Determination and Control\/}, Vol.~73,
  Springer, 2012.

\bibitem{KS1965}
Kustaanheimo, P., Schinzel, A., Davenport, H., and Stiefel, E.,
  \enquote{Perturbation Theory of Kepler Motion Based on Spinor
  Regularization,} {\em Journal f{\"u}r die Reine und Angewandte Mathematik\/},
  Vol.~1965, No. 218, apr 1965,
  pp.~204--219.~~\url{https://doi.org/10.1515/crll.1965.218.204}.

\bibitem{Chelnokov1981}
Chelnokov, Y.~N., \enquote{Regularization of the Equations of the Three
  Dimensional Two-Body Problem.} {\em Mechanics of Solids\/}, Vol.~16, No.~6,
  1981, pp.~10--18.

\bibitem{Vivarelli1983}
Vivarelli, M.~D., \enquote{The KS-Transformation in Hypercomplex Form,} {\em
  Celestial Mechanics and Dynamical Astronomy\/}, Vol.~29, No.~1, jan 1983,
  pp.~45--50.~~\url{https://doi.org/10.1007/bf01358597}.

\bibitem{Deprit1994}
Deprit, A., Elipe, A., and Ferrer, S., \enquote{Linearization: Laplace vs.
  Stiefel,} {\em Celestial Mechanics and Dynamical Astronomy\/}, Vol.~58,
  No.~2, feb 1994, pp.~151--201.~~\url{https://doi.org/10.1007/bf00695790}.

\bibitem{Vrbik1994}
Vrbik, J., \enquote{Celestial Mechanics via Quaternions,} {\em Canadian Journal
  of Physics\/}, Vol.~72, No. 3-4, mar 1994,
  pp.~141--146.~~\url{https://doi.org/10.1139/p94-023}.

\bibitem{Vrbik1995}
Vrbik, J., \enquote{Perturbed Kepler Problem in Quaternionic Form,} {\em
  Journal of Physics A: Mathematical and General\/}, Vol.~28, No.~21, nov 1995,
  pp.~6245--6252.~~\url{https://doi.org/10.1088/0305-4470/28/21/027}.

\bibitem{Waldvogel2006}
Waldvogel, J., \enquote{Quaternions and the Perturbed Kepler Problem,} {\em
  Celestial Mechanics and Dynamical Astronomy\/}, Vol.~95, No. 1-4, aug 2006,
  pp.~201--212.~~\url{https://doi.org/10.1007/s10569-005-5663-7}.

\bibitem{Waldvogel2008}
Waldvogel, J., \enquote{Quaternions for Regularizing Celestial Mechanics: the
  Right Way,} {\em Celestial Mechanics and Dynamical Astronomy\/}, Vol.~102,
  No. 1-3, mar 2008,
  pp.~149--162.~~\url{https://doi.org/10.1007/s10569-008-9124-y}.

\bibitem{Saha2009}
Saha, P., \enquote{Interpreting the Kustaanheimo--Stiefel Transform in
  Gravitational Dynamics,} {\em Monthly Notices of the Royal Astronomical
  Society\/}, Vol.~400, No.~1, nov 2009,
  pp.~228--231.~~\url{https://doi.org/10.1111/j.1365-2966.2009.15437.x}.

\bibitem{Broucke1971}
Broucke, R., Lass, H., and Ananda, M., \enquote{Redundant Variables in
  Celestial Mechanics,} {\em Astronomy and Astrophysics\/}, Vol.~13, aug 1971,
  pp.~390--398.

\bibitem{Deprit1975}
Deprit, A., \enquote{Ideal Elements for Perturbed Keplerian Motions,} {\em
  Journal of Research of the National Bureau of Standards, Section B:
  Mathematical Sciences\/}, Vol.~79B, No.~1, jan 1975,
  pp.~1--15.~~\url{https://doi.org/10.6028/jres.079b.001}.

\bibitem{Gurfil2005}
Gurfil, P., \enquote{Euler Parameters as Nonsingular Orbital Elements in
  Near-Equatorial Orbits,} {\em Journal of Guidance, Control, and Dynamics\/},
  Vol.~28, No.~5, sep 2005,
  pp.~1079--1084.~~\url{https://doi.org/10.2514/1.14760}.

\bibitem{Chelnokov2001}
Chelnokov, Y.~N., \enquote{The Use of Quaternions in the Optimal Control
  Problems of Motion of the Center of Mass of a Spacecraft in a Newtonian
  Gravitational Field: I,} {\em Cosmic Research\/}, Vol.~39, No.~5, 2001,
  pp.~470--484.~~\url{https://doi.org/10.1023/a:1012345213745}.

\bibitem{Chelnokov2003}
Chelnokov, Y.~N., \enquote{The Use of Quaternions in the Optimal Control
  Problems of Motion of the Center of Mass of a Spacecraft in a Newtonian
  Gravitational Field: II,} {\em Cosmic Research\/}, Vol.~41, No.~1, 2003,
  pp.~85--99.~~\url{https://doi.org/10.1023/a:1022359831200}.

\bibitem{Chelnokov2013}
Chelnokov, Y.~N., \enquote{Quaternion Regularization in Celestial Mechanics and
  Astrodynamics and Trajectory Motion Control. I,} {\em Cosmic Research\/},
  Vol.~51, No.~5, sep 2013,
  pp.~350--361.~~\url{https://doi.org/10.1134/s001095251305002x}.

\bibitem{Chelnokov2014}
Chelnokov, Y.~N., \enquote{Quaternion Regularization and Trajectory Motion
  Control in Celestial Mechanics and Astrodynamics: II,} {\em Cosmic
  Research\/}, Vol.~52, No.~4, jul 2014,
  pp.~304--317.~~\url{https:doi.org/10.1134/s0010952514030022}.

\bibitem{Chelnokov2019}
Chelnokov, Y.~N., \enquote{Perturbed Spatial Two-Body Problem: Regular
  Quaternion Equations of Relative Motion,} {\em Mechanics of Solids\/},
  Vol.~54, No.~2, mar 2019,
  pp.~169--178.~~\url{https://doi.org/10.3103/s0025654419030075}.

\bibitem{Libraro2014}
Libraro, P., Kasdin, N.~J., Choueiri, E.~Y., and Dutta, A.,
  \enquote{Quaternion-Based Coordinates for Nonsingular Modeling of
  High-Inclination Orbital Transfer,} {\em Journal of Guidance, Control, and
  Dynamics\/}, Vol.~37, No.~5, sep 2014,
  pp.~1638--1644.~~\url{https://doi.org/10.2514/1.g000613}.

\bibitem{Roa2017}
Roa, J. and Kasdin, N.~J., \enquote{Alternative Set of Nonsingular Quaternionic
  Orbital Elements,} {\em Journal of Guidance, Control, and Dynamics\/},
  Vol.~40, No.~11, nov 2017,
  pp.~2737--2751.~~\url{https://doi.org/10.2514/1.g002753}.

\bibitem{Bate1}
Bate, R.~R., Mueller, D.~D., and White, J.~E., {\em Fundamentals of
  Astrodynamics\/}, Dover Publications, 1971.

\bibitem{Clarke2}
Rao, A. and Clarke, K., \enquote{Performance optimization of a maneuvering
  re-entry vehicle using a legendre pseudospectral method,} {\em AIAA
  Atmospheric Flight Mechanics Conference and Exhibit\/}, American Institute of
  Aeronautics and Astronautics ({AIAA}), aug
  2002.~~\url{https://doi.org/10.2514/6.2002-4885}.

\bibitem{Clarke1}
Clarke, K., {\em Performance Optimization of a Common Aero Vehicle Using a
  Legendre Pseudospectral Method\/}, Master's thesis, Department of Aeronautics
  and Astronautics, Massachusetts Institute of Technology, Cambridge,
  Massachusetts, 2003.

\bibitem{Patterson2014}
Patterson, M.~A. and Rao, A.~V., \enquote{$\mathbb{GPOPS-II}$, {A} {MATLAB}
  {S}oftware for {S}olving {M}ultiple-{P}hase {O}ptimal {C}ontrol {P}roblems
  {U}sing $hp$-{A}daptive {G}aussian {Q}uadrature {C}ollocation {M}ethods and
  {S}parse {N}onlinear {P}rogramming,} {\em ACM Transactions on Mathematical
  Software\/}, Vol.~41, No.~1, oct 2014,
  pp.~1--37.~~\url{https://doi.org/10.1145/2558904}.

\bibitem{netlib}
Netlib, \enquote{fdlibm,} \url{https://www.netlib.org/fdlibm}, 1995, Accessed:
  2020-05-26.

\end{thebibliography}

\renewcommand{\baselinestretch}{1.5}
\normalsize\normalfont

\end{document}